\documentclass[12pt]{amsart}
\usepackage{mathrsfs, amssymb, array}
\usepackage[colorlinks]{hyperref}
\usepackage{mathcomp}
\usepackage[all]{xy}
\usepackage{longtable}
 
\newcommand{\comment}[1]{}

\newcommand{\EEE}{{\mathscr{E}}}
\newcommand{\OOO}{\mathscr{O}}

\newcommand{\FFF}{\mathcal{F}}

\newcommand{\CC}{\mathbb{C}}
\newcommand{\ZZ}{\mathbb{Z}}
\newcommand{\PP}{\mathbb{P}}
\newcommand{\QQ}{\mathbb{Q}}
\newcommand{\FF}{\mathbb{F}}
\newcommand{\RR}{\mathbb{R}}

\newcommand{\e}{\mathbf{e}}
\newcommand{\h}{\mathbf{h}}

\newcommand{\Sym}{\mathfrak{S}}

\newcommand{\xref}[1]{{\rm \ref{#1}}}

\newcommand{\comp}{\mathbin{\scriptstyle{\circ}}}

\newcommand{\Sing}{\operatorname{Sing}}

\newcommand{\Pic}{\operatorname{Pic}}
\newcommand{\Cl}{\operatorname{Cl}}
\newcommand{\pt}{\operatorname{pt}}
\newcommand{\Aut}{\operatorname{Aut}}
\newcommand{\Cr}{\operatorname{Cr}}

\newcommand{\rk}{\operatorname{rk}}

\newcommand{\Gr}{\operatorname{Gr}}

\newcommand{\dd}{\operatorname{d}}
\newcommand{\rr}{\operatorname{r}}

\newcommand{\pp}{\operatorname{p}}

\newcommand{\s}{\operatorname{s}}
\newcommand{\hh}{\operatorname{h}}
\newcommand{\W}{\operatorname{W}}

\newcommand{\Eu}{\operatorname{Eu}}

\renewcommand{\xref}[1]{\textup{\ref{#1}}}

\renewcommand\labelenumi{(\roman{enumi})}

\newtheorem{theorem}{Theorem}
\renewcommand{\thetheorem}{\thesection.\arabic{theorem}}
\numberwithin{theorem}{section}
\numberwithin{equation}{theorem}
\renewcommand{\theequation}{\thetheorem.\arabic{equation}}

\newtheorem{mtheorem}[theorem]{}
\newtheorem{stheorem}[theorem]{}

\theoremstyle{definition}

\newtheorem{case}[theorem]{}
\newtheorem{scase}[theorem]{}
\newtheorem*{remark*}{Remark}

\newcounter{NO}\numberwithin{NO}{theorem}

\def\numer{\refstepcounter{equation}{\rm (\theequation)}}

\newcounter{NOa}[NO]\numberwithin{NOa}{equation}

\newcounter{NN}\numberwithin{NN}{theorem}
\renewcommand{\theNN}{\arabic{NN}${}^o$}
\def\nr{\refstepcounter{NN}{\theNN}}%

\title{$G$-Fano threefolds, I}

\author{Yuri Prokhorov}
\thanks{The work was partially supported by RFBR grants Nos. 11-01-00336-a and
11-01-92613-KO\_a,
Leading Scientific Schools, No. 4713.2010.1,
and 
AG Laboratory HSE, RF government grant ag. 11.G34.31.0023}
\address{Department 
of Algebra, Faculty of Mathematics, Moscow State
University, Moscow, 119 991, RUSSIA
\newline\indent
Laboratory of Algebraic Geometry, SU-HSE, 
7 Vavilova Str., Moscow, 117312, RUSSIA
}
\email{prokhoro@gmail.com}

\begin{document}

\begin{abstract}
We classify Fano threefolds with only terminal 
singularities whose 
canonical class is Cartier and divisible by $2$, and 
satisfying an additional assumption that the 
$G$-invariant part of the Weil divisor class group is of rank $1$
with respect to an action of 
some group $G$.
In particular, we find a lot of examples
of Fano $3$-folds with ``many'' symmetries.
\end{abstract}
\maketitle
\tableofcontents

\section{Introduction.}

\begin{case}
Let $Y$ be an algebraic variety $X$ over a 
field $\Bbbk$ and let $G$ be a group.
Following works of Yu. I. Manin \cite{Manin-1967}
we say that $X$ is a $G$-variety if 
the group $G$ acts on $\bar X:=X\otimes_\Bbbk \bar \Bbbk$, where 
$\bar \Bbbk$ is the algebraic closure of $\Bbbk$.
Moreover, we assume that $X$, $G$ and $\Bbbk$ satisfy
one of the following two conditions.

(a) \textit{Algebraic case.} 
$G$ is the Galois group of $\bar \Bbbk$ over $\Bbbk$
acting on $\bar X$ through the second factor.
The action of $G$ on $X$ is trivial.

(b) \textit{Geometric case.}
The field $\Bbbk$ is algebraically closed,
$G$ is a finite group, 
and the action of $G$ on $X$ is given by a 
homomorphism $G\to \Aut_{\Bbbk}(X)$.

A $G$-morphism (resp. rational $G$-map) of 
$G$-varieties 
is a $\Bbbk$-morphism (resp. $\Bbbk$-rational map) 
commuting with the action of $G$ in the geometric case.
A projective $G$-morphism $f: X\to Z$ of normal 
$G$-varieties is called \textit{$G$-Mori fiber space }if 
$X$ has at worst 
terminal $G\QQ$-factorial singularities (see \ref{Terminal-singularities}), $f_*\OOO_X=\OOO_Z$,
the relative invariant Picard group $\Pic(X/Z)^G$ is of rank $1$, and 
the anticanonical divisor $-K_X$ is $f$-ample. 
In the particular case where $\dim Z=0$, $X$ is called a
\textit{$G\QQ$-Fano variety}. 
If furthermore the canonical divisor is Cartier, then we say 
that $X$ is $G$-Fano variety.

Throughout this paper we assume that the ground field 
has characteristic $0$. The following is an easy consequence 
of the Minimal Model Program \cite[0.3.14]{Mori-1988} (cf. \cite[4.2]{Prokhorov2009e}).

\begin{mtheorem}{\bf Proposition.}
Let $V$ be a $G$-variety of dimension $\le 3$.
The following are equivalent:
\begin{enumerate}
 \item 
$\kappa(V)=-\infty$,
 \item 
$V$ is geometrically uniruled,
 \item
$V$ is $G$-birationally isomorphic to a
variety $X$ having a structure of $G$-Mori fiber space.
\end{enumerate}
\end{mtheorem}

Birational classification of $G$-surfaces is developed very well
\cite{Manin-1967}, \cite{Iskovskikh-1979s-e}.
In this and subsequent papers we consider 
$G$-Mori fiber spaces $X\to Z$, where $\dim X=3$ and $Z$ is a point, i.e.
the case of $G\QQ$-Fano threefolds.
\end{case}

\begin{case}
\label{notation}
Let $X$ be a $G$-Fano threefold. 
It is well-known that $\Pic(X)$ is a finitely generated 
torsion free abelian group (see, e.g. \cite[Prop. 2.1.2]{Iskovskikh-Prokhorov-1999}).
Consider the following composed object:
\[
\operatorname{V}(X)=\bigl(\Cl(X),\Pic(X),K_X,(\phantom{v},\phantom{v},\phantom{v})\bigr), 
\]
where $\Pic(X)$ is regarded as a sublattice of $\Cl(X)$, 
$K_X\in \Pic(X)$ is the canonical class of $X$, and $(\phantom{v},\phantom{v},\phantom{v})$ is
the intersection form $\Pic(X)\times \Pic(X)\times\Cl(X)\to \ZZ$.
Since the singularities of $X$ are isolated cDV \cite{Reid-YPG1987}, $\Pic(X)$ is a primitive sublattice in 
$\Cl(X)$, i.e. the quotient $\Cl(X)/\Pic(X)$ is torsion 
free \cite[5.1]{Kawamata-1988-crep}. Moreover, since $\rho(X)^G=1$, we have
\begin{equation}
\text{$\Cl(X)^G$ is a subgroup of rank $1$
containing $K_X$.}
\end{equation}
\end{case}

 \begin{case}
In this paper we give a classification of one class 
of Gorenstein $G$-Fano threefolds. More precisely, we consider Fano threefolds
such that $-K_X=2S$ for some ample 
Cartier divisor $S$. Then $X$ is called a \textit{del Pezzo threefold}
(see \ref{definition-del-Pezzo}).
Smooth del Pezzo threefolds were classified by Iskovskikh
\cite{Iskovskikh-1980-Anticanonical}, see also \cite{Fujita-all},
\cite{Iskovskikh-Prokhorov-1999}. Singular ones were
discussed from different points of view in many works
\cite{Fujita-1986-1}, \cite{Shin1989}, \cite{Fujita1990}, \cite{Casagrande2008}, \cite{Jahnke2008}.
We are interested basically in group actions on terminal del Pezzo threefolds $X$
and the structure of the lattice $\Cl(X)$.
\end{case}

\begin{case}
\label{delta}
Let $S$ be a smooth del Pezzo surface of degree $d:=K_S^2$.
Then we have $\Pic(S) = \ZZ^{10-d}$. 
Define
\[
\Delta:=\{\alpha \in \Pic(S) \mid \alpha^2=-2,\quad \alpha\cdot K_S=0\}. 
\]
Then $\Delta$ is a root system in $(K_S)^{\perp}\otimes \RR$. Depending on $d$,
$\Delta$ is of the following type (\cite{Manin-Cubic-forms-e-II}): 
\bigskip
\begin{center}
\begin{tabular}[]{c|ccccccccc}
\hline
\\[-8pt]
d &1&2&3&4&5&6&7&$8'$&$8''$
\\
\hline
\\[-5pt]
$\Delta$ & $E_8$& $E_7$& $E_6$& $D_5$& $A_4$& $A_1\times A_2$&$-$&$A_1$&$-$
\\
\hline
\end{tabular}
\end{center}
where $8'$ (resp. $8''$) corresponds to $\PP^1\times\PP^1$
(resp. Hirzebruch surface $\FF_1$).
\end{case}

\begin{case}
\label{delta-1}
Now let $X$ be a del Pezzo threefold.
Let $S\in |{-}\frac12K_X|$ be a smooth member  \cite{Shin1989}
and let $\iota: S \hookrightarrow X$ be the natural embedding.
Then $S$ is a del Pezzo surface of degree $d=-\frac 18 K_X^3$.
It is easy to show that the restriction map 
$\iota^*: \Cl(X)\to \Pic(S)$ is injective and its cokernel
is torsion free (see
\ref{corollary-embedding}).
Define the following subsets in $\Delta\subset\Pic(S)$:
\[
\begin{array}{lll}
\Delta'&:=&\{\alpha \in \Pic(S) \mid\alpha^2=-2,\ \alpha\cdot K_S= \alpha\cdot \iota^*\Cl(X)=0\}, 
\\[5pt]
\Delta''&:=&\{\alpha \in \iota^*\Cl(X) \mid\alpha^2=-2,\ \alpha\cdot K_S=0 \}. 
\end{array}
\]
In other words, 
\[
\Delta'=\Delta\cap (\iota^*\Cl(X))^{\perp}, \qquad 
 \Delta''=\Delta\cap \iota^*\Cl(X).
\]
If $\Delta'$ (resp. $\Delta''$) is non-empty, then it is a root subsystem in $\Delta$.
Assume that $X$ is a $G$-variety. Then the group $G$ naturally acts on $\iota^*\Cl(X)$
and $\Delta''$ preserving the class of $K_S$
and the intersection form.
Our classification of $G$-del Pezzo threefolds is by types of 
root systems $\Delta'$ and $\Delta''$.
\end{case}

\begin{mtheorem} {\bf Theorem.}\label{main-theorem}
Let $X$ be a $G$-del Pezzo threefold and let $\dd(X):=-\frac18K_X^3$. 
There are the following possibilities:
\end{mtheorem}

\setlongtables \renewcommand{\arraystretch}{1.2}
\begin{longtable}{c|c|c|c|c|c|c|c|c}
\hline
&$\rr$ & $X$ & $\bar X$&$Z$ & $\Delta'$ & $\Delta''$ & $\pp$ & $\s$\tabularnewline
\endfirsthead
&$\rr$ & $X$ & $\bar X$ &$Z$& $\Delta'$ & $\Delta''$ & $\pp$ & $\s$ \tabularnewline
\hline
\endhead

\hline 
\multicolumn{9}{c}{$\dd(X)=1$} \tabularnewline[4pt]
\hline 
\nr\label{class-d=1-factorial} 
&$1$ & $V_{1}$ & $-$ & $\pt$ & $E_8$ & $-$ & $0$ & $21-\hh$ \tabularnewline[4pt]
\nr &$2$ & \eqref{rkCl=2-deg1-quadric}&$-$& $\PP^{1}$ & $D_7$ & $-$ & $0$ & $22-\hh$ \tabularnewline
\nr &$2$ & \eqref{rkCl=2-P1-bundle-degree=1}&$-$& $\PP^{2}$ & $A_7$ & $-$ & $0$ &$22$ \tabularnewline

\nr\label{class-d=1-r=2-birational} &$2$ &\eqref{rkCl=2-deg1-birational}& $V_{2}$ & $\pt$ & $E_7$ & $A_{1}$ & $2$ &$22-\hh$, $\hh\le10$ \tabularnewline

\nr &$3$ & &\eqref{rkCl=2-deg2-quadric}& $\PP^{1}$ & $D_6$ & $2A_1$ & $4$ & $23-\hh$ \tabularnewline

\nr &$3$ & &\eqref{rkCl=2-P1-bundle-degree=2}& $\PP^{2}$ & $A_6$ & $A_{1}$ & $2$ & $23$\tabularnewline

\nr\label{rkCl=3-V3} &$3$ && $V_{3}$ & $\pt$ & $E_6$ & $A_2$ & $6$ & $23-\hh$ \tabularnewline

\nr &$4$ & &\eqref{rkCl=2-P1-bundle-degree=3}& $\PP^{2}$ & $A_5$ & $A_1\times A_2$ & $8$ & $24$\tabularnewline

\nr &$4$ && $V_{4}$ & $\pt$ & $D_5$ & $A_3$ & $12$ &$24-\hh$, $\hh\le 2$ \tabularnewline

\nr &$5$ & &\eqref{rkCl=2-deg4-quadric}& $\PP^{1}$ & $D_4$ & $D_4$ & $24$ & $25-\hh$ \tabularnewline

\nr &$5$ && $V_{5}$ & $\pt$ & $A_4$ & $A_4$ & $20$ & $25$\tabularnewline

\nr &$6$ & &\eqref{rkCl=2-degree=5}& $\PP^{2}$ & $A_3$ & $D_5$ & $40$ & $26$ \tabularnewline

\nr &$7$ & \xref{main-theorem-semimaximal} &$V_{6}$ & $\PP^{2}$ & $A_2$ & $E_6$ & $72$ & $27$\tabularnewline

\nr &$8$ &\xref{maximal-rk-d=1} &$\PP^{3}$ & $\pt$ & $A_{1}$ & $E_{7}$ & $126$ & $28$\tabularnewline

\hline 
\multicolumn{9}{c}{$\dd(X)=2$} \tabularnewline
\hline 
\nr\label{class-d=2-factorial} &$1$ & $V_{2}$ &$-$ & $\pt$ & $E_7$ & $-$ & $0$ & $10-\hh$\tabularnewline

\nr &$2$ & \eqref{rkCl=2-deg2-quadric} &$-$& $\PP^{1}$ & $D_6$ & $A_{1}$ & $0$ & $11-\hh$ \tabularnewline

\nr &$2$ & \eqref{rkCl=2-P1-bundle-degree=2} &$-$& $\PP^{2}$ & $A_6$ & $-$ & $0$ & $11$\tabularnewline

\nr\label{class-d=2-r=2-birational} &$2$ &\eqref{rkCl=2-deg2-birational}& $V_{3}$ & $\pt$ & $E_6$ & $-$ & $2$ & $11-\hh$, $\hh\le 5$ \tabularnewline

\nr &$3$ & \xref{theorem-primitive-3}\xref{primitive-rho=3-d=2} &$-$& $(\PP^1)^2$ & $A_5$ & $A_2$ & $0$ & $12$\tabularnewline

\nr &$3$ & &\eqref{rkCl=2-P1-bundle-degree=3}& $\PP^{2}$ & $A_5$ & $A_{1}$ & $2$ & $12$\tabularnewline

\nr &$3$ & &$V_{4}$ & $\pt$ & $D_5$ & $A_{1}$ & $4$ & $12-\hh$, $\hh\le 2$\tabularnewline

\nr\label{class-2-d4} &$4$ & &\eqref{rkCl=2-deg4-quadric} &$\PP^{1}$ & $D_4$ & $3A_1$ & $8$ & $13-\hh$ \tabularnewline

\nr &$4$ && $V_{5}$ & $\pt$ & $A_4$ & $A_2$ & $6$ & $13$ \tabularnewline

\nr &$5$ && \eqref{rkCl=2-degree=5} & $\PP^{2}$ & $A_3$ & $A_1\times A_3$ & $12$ & $14$\tabularnewline

\nr &$6$ &\xref{main-theorem-semimaximal} & $V_{6}$& $\PP^{2}$ & $A_2$ & $A_5$ & $20$ & $15$\tabularnewline

\nr\label{class-maximal-rk-d=2} &$7$ &\xref{maximal-rk-d=2}& $\PP^{3}$ & $\pt$ & $A_{1}$ & $D_{6}$ & $32$ & $16$\tabularnewline

\hline
\multicolumn{9}{c}{$\dd(X)=3$} \tabularnewline
\hline 
\nr\label{class-cubic-factorial} &$1$ & $V_{3}$ &$-$ & $\pt$ & $E_6$ & $-$ & $0$ & $\le 5$\tabularnewline

\nr\label{class-cubic-P1-bundle} &$2$ & \eqref{rkCl=2-P1-bundle-degree=3} &$-$& $\PP^{2}$ & $A_5$ & $A_{1}$ & $0$ & $6$ \tabularnewline

\nr\label{class-cubic-empty-delta} &$3$ & \xref{construction-d=3-empty-delta} &\eqref{rkCl=2-deg4-quadric}& $\PP^{1}$ & $D_4$ & $-$ & $3$ & $6\le \s\le 7$ \tabularnewline

\nr\label{class-cubic-subSegre} &$5$ &\xref{main-theorem-semimaximal} & $V_{6}$& $\PP^{2}$ & $A_2$ & $2A_2$ & $9$ & $9$\tabularnewline

\nr\label{class-cubic-Segre} &$6$ &\xref{Segre-cubic}& $\PP^{3}$ & $\pt$ & $A_{1}$ & $A_{5}$ & $15$ & $10$ \tabularnewline

\hline
\multicolumn{9}{c}{$\dd(X)=4$} \tabularnewline
\hline 
\nr &$1$ & $V_{4}$ &$-$ & $\pt$ & $D_5$ & $-$ & $0$ & $\le 2$ \tabularnewline

\nr &$2$ & \eqref{rkCl=2-deg4-quadric} &$-$& $\PP^{1}$ & $D_4$ & $-$ & $0$ & $1\le \s\le 3$ \tabularnewline

\nr &$3$ & \xref{theorem-primitive-3}\xref{primitive-rho=3-d=4}&$-$& $(\PP^1)^2$ & $A_3$ & $2A_1$ & $0$ & $4$\tabularnewline

\nr &$4$ &\xref{main-theorem-semimaximal} & $V_{6}$& $\PP^{2}$ & $A_2$ & $2A_1$ & $4$ & $5$\tabularnewline

\nr &$5$ & \xref{maximal-deg=4}&$\PP^{3}$ & $\pt$ & $A_{1}$ & $A_{1}\times A_{3}$ & $8$ & $6$ \tabularnewline

\hline
\multicolumn{9}{c}{$\dd(X)=5$} \tabularnewline
\hline 
\nr &$1$ & $V_{5}$&$-$ & $\pt$ & $A_4$ & $-$ & $0$ & $0$ \tabularnewline

\hline 
\multicolumn{9}{c}{$\dd(X)=6$} \tabularnewline
\hline 
\nr &$2$ & $V_{6}$ &$-$& $\PP^{2}$ & $A_2$ & $A_{1}$ & $0$ &$0$ \tabularnewline

\nr &$3$ & $(\PP^{1})^{3}$ &$-$ & $\pt$ & $A_{1}$ & $A_2$ & $0$ & $0$ \tabularnewline

\hline
\multicolumn{9}{c}{$\dd(X)=8$} \tabularnewline
\hline 
\nr &$1$ & $\PP^{3}$&$-$ & $\pt$ & $-$ & $-$ & $0$ &$0$ \tabularnewline
\hline
\end{longtable}
Here $\bar X/Z$ is a primitive birational 
model of $X$ (see Theorem \ref{structure-theorem-imprimitive})
and $\hh:=h^{1,2}(\hat X)$, where $\hat X$ is the standard resolution of $X$.
For compactness, we denote $(\PP^1)^k:=\underbrace{\PP^1\times \cdots \times \PP^1}_k$.
For other notation we refer to \ref{Notation}.

\begin{case}{\bf Remark.}
For $\dd(X)\le 2$ any del Pezzo threefold automatically has $G$-structure (see Remark \ref{def-Bertini-Geiser}).
So, in this case, \ref{class-d=1-factorial} -- \ref{class-maximal-rk-d=2}
is a complete list of del Pezzo threefolds with $\dd(X)\le 2$.
\end{case}

\begin{case}{\bf Remark.}
Singular three-dimensional cubics (without group action) 
whose singularities are only \textit{nodes} and their small resolutions were classified in \cite{Finkelnberg}.
There is the following correspondence between our list and the classification in \cite{Finkelnberg}:
\ref{class-cubic-Segre} $\longleftrightarrow$ J15,
\ref{class-cubic-subSegre} $\longleftrightarrow$ J14,
\ref{class-cubic-empty-delta} $\longleftrightarrow$ J11,
\ref{class-cubic-P1-bundle}$\longleftrightarrow$ J9,
\ref{class-cubic-factorial}$\longleftrightarrow$ J1--J5.
\end{case}

We hope that our result can be useful for applications to 
the classification of finite subgroups of the Cremona group
$\Cr_3(\Bbbk)$ \cite{Prokhorov2009e}, \cite{Prokhorov2009d}, and also the the birational classification of 
rational algebraic threefolds over non-closed fields 
(cf. \cite{Manin-1967}).


\noindent
\textbf{Acknowledgments.}
This paper is an extended version of my talk given at 
``Workshop on Moduli and Birational Geometry''
at POSTECH, Korea in August, 2010.
I would like to thank POSTECH for the support and hospitality.
I also would like to acknowledge discussions that I have on this subject 
with I. Arzhantsev, I. Cheltsov, A. Kuznetsov, K. Shramov
and D. Timashev.
Thanks are also due to the 
referee for suggesting several improvements.

\section{Preliminaries.}
\label{section-Preliminaries}

\begin{case}{\textbf{Notation.}}
\label{Notation}
We work over an algebraically closed field of 
characteristic $0$.
Throughout this paper $X$ denotes a del Pezzo threefold with at worst 
terminal Gorenstein
singularities.
Thus we can write $-K_X=2S$, where $S=S_X$ 
is a ample Cartier divisor of $S$ defined up to linear equivalence.
Everywhere below we use the following notation:
\begin{itemize}
\item[] $\rho=\rho(X):=\rk \Pic(X)$;
\item[] $\Cl(X)$ is the Weil divisor class group;
\item[] $\rr=\rr(X):=\rk \Cl(X)$;
\item[] $\dd=\dd(X):=S^3=-K_X^3/8$, the degree of $X$;
\item[] $\pp=\pp(X)$ is the number of planes on $X$;
\item[] $\s=\s(X)$ is the number of singular points of $X$ under an additional assumption that 
$X$ has at worst nodes;
\item[] $V_6\subset \PP^7$ is a \textit{smooth} del Pezzo threefold 
with $\dd(X)=6$ and $\rho=2$, see Theorem \ref{th-del-Pezzo-2}; 
\item[] $V_5\subset \PP^6$ is a \textit{smooth} del Pezzo threefold of degree $5$
(see \cite{Iskovskikh-Prokhorov-1999}); 
\item[] $V_d$, for $d=1,\dots, 4$, is a del Pezzo threefold of degree $d$
with terminal \textit{factorial} singularities (see Theorem \ref{th-del-Pezzo-2}).
\end{itemize}
\end{case}
\begin{case}
\label{Terminal-singularities}
\textbf{Terminal singularities (see \cite{Reid-YPG1987}).} 
Let $(X, P)$ be a germ of a three-dimensional terminal singularity.
Then $(X, P)$ is isolated, i.e, $\Sing(X)=\{P\}$.
If additionally $(X, P)$ is Gorenstein, then
$(X, P)$ is analytically isomorphic to a hypersurface 
singularity of multiplicity $2$.

Let $X$ be a threefold with Gorenstein terminal singularities.
Then any Weil $\QQ$-Cartier divisor is Cartier (see e.g. \cite[Lemma 5.1]{Kawamata-1988-crep}). 
Equivalently, $\Pic(X)$ is a primitive sublattice in $\Cl(X)$. 

Let $X$ be a $G$-variety. We say that $X$ has only
\textit{$G\QQ$-factorial singularities} if any $G$-invariant 
Weil divisor is $\QQ$-Cartier.

\begin{stheorem}{\bf Theorem-Definition (\cite[Corollary 4.5]{Kawamata-1988-crep}).}
Let $X$ be a threefold with terminal singularities.
Then there exists a projective birational morphism
$\xi: \hat X\to X$ such that 
\begin{enumerate}
\item 
$\hat X$ is normal and has only terminal $\QQ$-factorial singularities;
\item 
$\xi$ is a crepant morphism, that is, $K_{\hat{X}}=\xi^*K_X$;
\item 
$\xi$ is small, that is, its exceptional locus does not contain any
divisors.
\end{enumerate}
Such a morphism is called \emph{$\QQ$-factorialization} of $X$.
Any two $\QQ$-factorializations of $X$ are connected by a sequence of flops.
\end{stheorem}
\end{case}

\begin{stheorem}{\bf Theorem \cite{Cutkosky-1988}.}
\label{theorem-Cutkosky}
Let $X$ be a rationally connected threefold with terminal factorial singularities.
Assume that $-K_X=2S$ for some divisor $S$ and $\rho(X)>1$.
Let $f: X\to Z$ be an extremal $K_X$-negative contraction.
Then one of the following holds:
\begin{enumerate}
\item 
$Z\simeq \PP^1$ and $f$ is a quadric bundle, i.e. 
there is an embedding $X\hookrightarrow \PP(\EEE)$,
where $\EEE$ is a rank $4$ vector bundle on $Z$, so that
each fiber of $f$ is a quadric in the fiber of $\PP(\EEE)/Z$; 
\item
$X$ is smooth, $Z$ is a smooth rational surface, and $X=\PP(\EEE)$, 
where $\EEE$ is a rank $2$ vector bundle on $Z$;
\item
$Z$ is a threefold with terminal factorial singularities 
and $f$ is blowup of a smooth point on $Z$.
\end{enumerate}
\end{stheorem}

\section{Generalities on del Pezzo threefolds}
\label{section-Generalities}
\begin{case}{\bf Definition.}
\label{definition-del-Pezzo}
Let $X$ be a projective variety $X$ with at worst 
terminal Gorenstein singularities.\footnote{In papers \cite{Fujita1990}, \cite{Casagrande2008}
authors considered del Pezzo varieties whose
singularities are more general than terminal.}
We say that $X$ is a \textit{del Pezzo threefold}
(resp. \textit{weak del Pezzo threefold}) if
its anti-canonical class $-K_X$ is divisible by $2$ and is ample 
(resp. nef and big).
\end{case}
\noindent
Note that if $X$ is a Fano threefold  with at worst 
terminal Gorenstein singularities such that $-K_X$ is divisible 
by some positive integer $q$, then $q\le 4$.
Moreover, $q=4$ iff $X\simeq \PP^3$ and $q=3$ iff $X$ is a quadric in $\PP^4$
 (see e.g. \cite[Th. 3.1.14]{Iskovskikh-Prokhorov-1999}).
Thus, for a del Pezzo threefold  $X\not \simeq \PP^3$,
the divisor $-\frac12 K_X$
is a primitive element of the lattice $\Cl(X)$.

\begin{mtheorem}{\bf Theorem 
(\cite{Iskovskikh-1980-Anticanonical}, \cite{Shin1989}, \cite{Fujita-all}, \cite{Fujita-1986-1}, \cite{Fujita1990}).}
\label{th-del-Pezzo-2}
Let $X$ be a del Pezzo threefold and let $S=-\frac{1}{2}K_X$.
\begin{enumerate}
\item 
$\dim |S|=\dd(X)+1$ and $\dd(X)\le 8$.
\item 
The linear system $|S|$
is base point free \textup(resp. very ample\textup) for $\dd(X)\ge 2$ 
\textup(resp. $\dd(X)\ge 3$\textup).
If $\dd(X)\ge 4$, then the image of $X_{\dd(X)}\subset \PP^{\dd(X)+1}$
of $X$ under the embedding given by $|S|$ is an intersection of quadrics.
\item
If $\dd(X)=1$, then the linear system $|S|$ has a unique base point which is a smooth point of $X$.
In this case $|S|$ defines a rational map $X \dashrightarrow \PP^{2}$
whose general fiber is an elliptic curve. The variety $X$ is isomorphic to a 
hypersurface of degree $6$ in the weighted projective space $\PP(1^3,2,3)$.
\item
If $\dd(X)=2$, then $|S|$ defines a double cover $X\to \PP^3$ whose 
branch locus $B\subset \PP^3$ is a surface of degree $4$ with at worst isolated singularities.
The variety $X$ is isomorphic to a 
hypersurface of degree $4$ in $\PP(1^{4},2)$.
\item
If $\dd(X)=3$, then $X$ is isomorphic to a cubic in $\PP^{4}$.
\item
If $\dd(X)=4$, then $X$ is isomorphic to a complete intersection of two quadrics
in $\PP^{5}$.

\item
If $\dd(X)= 5$, then  $\rho(X)=1$ \cite{Namikawa-1997}.
A smooth variety of this type is unique up to isomorphism
and isomorphic to a  section of 
the Grassmann variety $\Gr(2,5)\subset \PP(\wedge^2 \CC^5)$ by a subspace of codimension $3$.
\item
If $\dd(X)=6$, then $\rho(X)=3$ or $2$ \cite{Namikawa-1997}.
In both cases a smooth variety of this type is unique up to isomorphism and 
either $X\simeq \PP^1\times \PP^1\times \PP^1$ or 
$X$ is isomorphic to a divisor of bidegree $(1,1)$ in $\PP^2\times \PP^2$.

\item
If $\dd(X)= 7$, then 
$X\simeq \PP (\OOO_{\PP^2}\oplus \OOO_{\PP^2}(1))$.
\item
If $\dd(X)= 8$, then 
$X\simeq \PP^3$. 
\end{enumerate}
\end{mtheorem}

The del Pezzo threefolds with $\dd(X)=1$ and $\dd(X)=2$ have their names: 
\textit{double Veronese cone} and \textit{quartic double solid}, respectively.

\begin{scase}{\bf Remark.}
\label{def-Bertini-Geiser}
Let $X$ be a del Pezzo threefold of degree $1$ (resp. $2$).
Then there is a double cover $\varphi: X\to \PP(1^{3},2)$
(resp. $\varphi: X\to \PP^3$). The corresponding natural Galois involution
$X\to X$ is called \textit{Bertini} (resp. \textit{Geiser}) involution.
Therefore, any del Pezzo threefold $X$ with $\dd(X)\le 2$ is a $G$-del Pezzo
(in the geometric sense).
\end{scase}

\begin{mtheorem} {\bf Corollary.}
\label{lemma-Q-factorial}
Let $X$ be a del Pezzo threefold. If 
$X$ is factorial and singular, then $\rho(X)=1$.
\end{mtheorem}

\begin{proof}[Outline of the proof]
By Theorem \ref{th-del-Pezzo-2} we have to check only the case $\dd(X)=6$.
Alternatively, one can use Theorem \ref{theorem-Cutkosky}
(see the arxiv version).
\end{proof}

\begin{mtheorem}{\bf Proposition (see \cite{Iskovskikh-Prokhorov-1999}).}
Let $X$ be a smooth del Pezzo threefold. Then the Hodge number $h^{1,2}(X)$ 
is given by the following table:
\begin{center}
\begin{tabular}{|l||l|l|l|l|l|l|l|l|}
\hline
&&&&&&&&
\\[-5pt]
$\dd(X)$&
$1$&$2$&$3$&$4$&$5$&$6$&$7$&$8$
\\[5pt]
\hline
&&&&&&&&
\\[-5pt]
$h^{1,2}(X)$&
$21$& $10$ & $5$ & $2$ & $0$& $0 $& $0 $& $ 0$
\\[5pt]
\hline
\end{tabular}
\end{center}
\end{mtheorem}

\begin{case} {\bf Definition.}
Let $X$ be a weak del Pezzo threefold and let $S=-\frac12K_X$.
An irreducible surface $\Pi\subset X$ is called a \textit{plane} 
if $S^2\cdot \Pi=1$ and, in case $\dd(X)=1$, the base point of $|S|$
does not lie on $\Pi$.

\begin{stheorem}{\bf Lemma.}
Let $X$ be a del Pezzo threefold.
If $\Pi\subset X$ is a plane, then $\Pi\simeq \PP^2$
and $\OOO_{\Pi}(S)=\OOO_{\PP^2}(1)$.
\end{stheorem}

\begin{proof}
The statement is obvious if $\dd(X)\ge 3$ because 
the divisor $S$ is very ample in this case.
If $\dd(X)=2$, then $|S|$ defines a double cover $\varphi:X\to \PP^3$
so that $\varphi(\Pi)$ is a projective plane on $\PP^3$.
Thus $\varphi|_{\Pi}:\pi \to \varphi(\Pi)\simeq \PP^2$ is a finite birational
morphism, so it is an isomorphism.
Finally if $\dd(X)=1$, then $|S|$ defines a rational map $\varphi:X\dashrightarrow \PP^2$
so that its restriction to $\Pi$ is a morphism which must be 
finite and birational. As above we get $\Pi\simeq \PP^2$.
\end{proof}

\begin{stheorem}{\bf Lemma.}\label{lemma-factorialization-plane}
If $\Pi\subset X$ is a plane, then there is a 
$\QQ$-factorialization $\xi:\hat X\to X$ such that
for the proper transform $\hat \Pi$ of $\Pi$ we have 
$\hat \Pi\simeq \PP^2$ and $\OOO_{\hat \Pi}(\hat \Pi)\simeq \OOO_{\PP^2}(-1)$.
Therefore, $\hat \Pi$ is contractible, i.e. there is a birational contraction 
$\hat X\to X'$ of $\hat \Pi$ to a smooth point.
Conversely, if $\xi:\hat X\to X$ is a 
$\QQ$-factorialization and 
$\hat \Pi\subset \hat X$ is an irreducible surface
such that $\hat \Pi\simeq \PP^2$ and $\OOO_{\hat \Pi}(\hat \Pi)\simeq \OOO_{\PP^2}(-1)$,
then $f(\hat \Pi)$ is a plane on $X$.
\end{stheorem}

\begin{proof} 
Let $\Pi\subset X$ be a plane.
Take a $\QQ$-factorialization $\xi:\hat X\to X$ so that $\hat \Pi$ is $f$-nef.
One can do it by performing flops over $X$.
Assume that $\hat \Pi$ is nef.
Then by the base point free theorem the linear system
$|n\hat \Pi|$ is base point free for $n\gg 0$.
Hence $|n\Pi|$ has no fixed components.
Since $X$ has at worst isolated singularities, by adjunction we have
\[
K_{\Pi}=(-2S+\Pi)|_{\Pi} \ge -2S|_{\Pi},
\]
a contradiction.

Thus $\hat \Pi$ is not nef. Then there is a $K_{\hat X}$-negative extremal ray 
$R$ such that $\hat \Pi\cdot R<0$.
Since $K_{\hat X}$ is divisible by $2$, from the classification 
of extremal rays (Theorem \ref{theorem-Cutkosky}) we see that
$\hat \Pi$ is contractible to a smooth point, 
$\hat \Pi\simeq \PP^2$ and $\OOO_{\hat \Pi}(\hat \Pi)\simeq \OOO_{\PP^2}(-1)$.
The converse statement is obvious.
\end{proof}

\begin{stheorem}{\bf Lemma.}\label{lemma-plane-restriction}
Let $X$ be a del Pezzo threefold and let $S\in |{-}\frac12 K_X|$ 
be a smooth member. Let $l\in \Pic (S)$ be an element such that
$l^2=l\cdot K_S=-1$ \textup(the class of a line $L\subset S$\textup).
Assume that $l\in \iota^*\Cl(X)$, where $\iota: S\hookrightarrow X$ is the embedding.
Then there exists a unique plane $\Pi\subset X$ such that 
$\iota^*\Pi=l$ \textup(i.e., $\Pi\cap S=L$\textup). 
\end{stheorem}
\begin{proof} 
Denote by $\Pi$ any divisor whose class coincides with $\iota^*l$.
Let $\xi:\hat X\to X$ be a $\QQ$-factorialization
as in Lemma \ref{lemma-factorialization-plane}, let $\hat S:=\xi^{-1}(S)$, 
and let $\hat \Pi$ be the proper transform of $\Pi$. 
By Shokurov's adjunction theorem the pair $(\hat X,\hat\Pi)$ is purely log terminal (PLT).
Hence, by the Kawamata-Viehweg vanishing \cite[Prop. 1]{Fukuda1997}
\[
H^1 (\hat X, \OOO_{\hat X}(\hat \Pi -\hat S)) =
H^1 (\hat X, \OOO_{\hat X}(\hat S+K_{\hat X}+\hat \Pi)) =0.
\]
Then one can see from the exact sequence
\[
0 \longrightarrow \OOO_{\hat X}(\hat \Pi -\hat S) \longrightarrow 
\OOO_{\hat X}(\hat \Pi) \longrightarrow \OOO_{\hat S}(\iota^*l) \longrightarrow 0
\]
that $H^0(\hat X, \OOO_{\hat X}(\hat \Pi))\neq 0$, so we may assume that
both $\hat \Pi$ and $\Pi$ are effective. Since $S^2\cdot \Pi=1$, $\Pi$ is a plane. 
Finally, if there is another plane $\Pi'$ such that 
$\iota^*\Pi'=l$, then $\Pi\sim \Pi'$ and $\OOO_{\hat\Pi}(\hat \Pi)=
\OOO_{\hat\Pi}(\hat \Pi')$ is positive, a contradiction.
\end{proof}

\begin{scase}{\bf Definition.}
We say that a del Pezzo threefold $X$ is \textit{imprimitive} if it contains 
at least one plane.
Otherwise we say that $X$ is \textit{primitive}.
\end{scase}
\end{case}

The following two theorems are easy consequences of 
{\cite[Prop. 2.8]{Casagrande2008}}.
\begin{mtheorem}{\bf Theorem.}
\label{structure-theorem-primitive}
Let $X$ be a primitive weak del Pezzo threefold with at worst terminal 
Gorenstein singularities. 
Let $\xi : \hat X\to X$ be a $\QQ$-factorialization.
Then there exists a $K_{\hat X}$-negative Mori contraction $f:\hat X\to Z$ 
such that one of the following holds:
\begin{enumerate}
\item 
$Z$ is a point, $\rho(\hat X)=1$, $X$ is factorial, and $\xi$ is an isomorphism;
\item 
$Z\simeq \PP^2$, $\rho(\hat X)=2$, and $f$ is a $\PP^1$-bundle,
i.e. $\hat X$ is smooth and $\hat X=\PP_{\PP^2}(\EEE)$, where $\EEE$ is a 
rank-$2$ vector bundle on $\PP^2$;
\item 
$Z\simeq \PP^1$, $\rho(\hat X)=2$, and $f$ is a 
quadric bundle,
i.e. any fiber of $f$ is an irreducible quadric in $\PP^3$;
\item 
$Z\simeq \PP^1\times \PP^1$, $\rho(\hat X)=3$, and $f$ is a $\PP^1$-bundle.
\end{enumerate}
\end{mtheorem}
\begin{proof}
Almost all the statements are proved in \cite[Prop. 2.8]{Casagrande2008}.
We have to show only that $Z\not \simeq \FF_2$.
Indeed, if $\dim Z=2$, then for a general member $\bar S\in |{-}\frac12 K_{\bar X}|$, 
the restriction $f|_{\bar S}: \bar S\to Z$
is birational. Hence $Z$ is a del Pezzo surface.
\end{proof}

\begin{mtheorem}{\bf Theorem.}
\label{structure-theorem-imprimitive}
Let $X$ be an imprimitive del Pezzo threefold with at worst terminal 
Gorenstein singularities. 
Then there exists a diagram
\begin{equation}
\label{eq-construction}
\xymatrix{
&
\hat{X} \ar[r]^{\sigma}\ar[ld]_{\xi} &
\bar{X} \ar[rd]^{f}
\\
X&&&Z
}
\end{equation}
where
\begin{enumerate} 
\item 
$\xi$ is a $\QQ$-factorialization;
 \item 
$\hat{X}$ is an weak del Pezzo threefold with at worst terminal 
factorial singularities;
 \item 
$\sigma$ is a blowup in smooth distinct points $P_1,\dots,P_n\in \bar X$;
 \item 
$\dd({X})=\dd(\hat{X})=\dd(\bar{X})+n$;

 \item \label{structure-theorem-imprimitive-v}
$\bar X$ is a primitive weak del Pezzo threefold with $\rho(\bar X)\le 2$,
thus $\bar X$ is described by \textup{(i)-(iii)} of Theorem \xref{structure-theorem-primitive}.
\end{enumerate}
\end{mtheorem}

\begin{stheorem}{\bf Corollary.}
\label{cor-2-rkcl}
Let $X$ be a del Pezzo threefold.
Then $\rr(X)+\dd(X)\le 9$.
\end{stheorem}
\begin{proof}
 We have
$
9\ge \rho(\bar X)+\dd(\bar X)=\rho(\hat X)+\dd(\hat X)=\rr(X)+\dd(X).
$
\end{proof}

\begin{stheorem}{\bf Corollary.}
\label{corollary-embedding}
Let $X$ be a weak del Pezzo threefold and let $S\in |{-}\frac12 K_X|$
be a smooth element.
Then the restriction map $\Cl(X)\to \Pic(S)$ is injective and its cokernel
is torsion free. 
\end{stheorem}
\begin{proof}
Clearly the assertion is invariant under taking 
small modifications.
In view of construction \eqref{eq-construction},
it is sufficient to prove that the restriction map 
$\Cl(\bar X)\to \Pic(\bar S)$ is injective and its cokernel
is torsion free, where $\bar S=\sigma(S)$.
Thus we may assume that $X$ is a primitive factorial weak del Pezzo
threefold.
The assertion is obvious if $\rho(X)=1$.
Assume that $\rho(X)=2$.
Then $\rho(Z)=1$. Let $\Theta$ be the ample generator of 
$\Pic(Z)$. The group
$\Cl(X)$ is generated by $f^* \Theta$ and the class of $S$.
Recall that $Z$ is either $\PP^1$ or $\PP^2$.
Hence $f^*\Theta|_S$ is either a conic or the pull-back of 
a line on $\PP^2$, respectively.
It is easy to see that $f^*\Theta|_S$ and $-K_S\sim S|_S$
generate a rank 2 primitive sublattice in $\Pic(S)$.
The case $\rho(X)=3$ can be treated similarly.
\end{proof}

\section{Primitive del Pezzo threefolds with $\rr(X)=3$}
\label{section-Primitive-3}
\begin{mtheorem}{\bf Lemma.}\label{lemma-primitive-r=3}
Let $X$ be a primitive del Pezzo threefold with $\rr(X)=3$
and let $\FFF=|F|$ be a complete one-dimensional linear system \textup(pencil\textup)
of Weil divisors without fixed components.
There is a small $\QQ$-factorialization $\xi: \hat X\to X$ such that 
the proper transform $\hat \FFF$ of $\FFF$ on $\hat X$ is 
base point free and defines a fibration $f: \hat X\to \PP^1$.
Moreover, $f$ factors through 
a \textup(not unique\textup) $\PP^1$-bundle contraction
\begin{equation}
\label{equation-primitive-r=3-fibration}
f: \hat X \overset{f_1}{\longrightarrow} \PP^1\times \PP^1 \overset{f_2}{\longrightarrow} \PP^1
\end{equation}
\end{mtheorem}

\begin{proof}
Take a $\QQ$-factorialization $\xi: \hat X\to X$ so that $\hat \FFF$ is
$\xi$-nef (one can get it by performing flops over $X$).
Then $\hat \FFF$ is nef. Indeed, otherwise there is a $K_{\hat X}$-negative extremal ray 
$R$ such that $\hat \FFF\cdot R<0$. Since $\hat \FFF$ has no fixed components,
$R$ must define a flipping contraction. 
On the other hand, $K_X$ is Cartier, a contradiction \cite[Th. 6.2]{Mori-1988}.
Thus $\hat \FFF$ is nef. Then $\hat \FFF$ defines a contraction to
a (rational) curve by the base point 
free theorem.
Further, since $\rr(X)=3$, we have $\rho(\hat X)=3$.
Running the MMP over $\PP^1$ we obtain $f_1$.
\end{proof}

\begin{scase}
{\bf Remark-definition.}
In notation of \eqref{equation-primitive-r=3-fibration}, 
another ruling on $\PP^1\times \PP^1$ defines another pencil $\FFF'$ on $X$. 
In this situation, we say that pencils $\FFF$ and $\FFF'$ are \textit{conjugate}.
Thus there is one-to-one correspondence between 
\begin{enumerate}
\item
the set of pairs of conjugate pencils $\FFF$, $\FFF'$
and 
\item
the set of $\QQ$-factorializations $X'\to X$ together with a
structure of $\PP^1$-bundle $f' : \hat{X'}\to \PP^1\times \PP^1$.
\end{enumerate}
\end{scase}

\begin{stheorem}{\bf Corollary.}
The cone of effective divisors 
$\overline{\operatorname{NE}}^1(X)$ is generated by 
classes of pencils $\FFF$ as in Lemma \xref{lemma-primitive-r=3}.
\end{stheorem}
\begin{proof}
Let $\xi :\hat{X}\to X$ be a small $\QQ$-factorialization.
There are natural identifications 
$\Cl(X)= \Cl(\hat{X})$ and $\overline{\operatorname{NE}}^1(X)=\overline{\operatorname{NE}}^1(\hat{X})$.
The variety $\hat{X}$ is a Mori dream space \cite{Hu2000}. Hence 
$\overline{\operatorname{NE}}^1(\hat{X})$ is a  polyhedral cone generated by a 
finite number of effective divisors $D_i$.
Running $D_i$-MMP on $\hat X$, after a number of flops, 
we get a $\PP^1$-bundle over 
$\PP^1\times \PP^1$ (because $X$ is primitive).
This shows that $D_i$ must coincide with some $\FFF$.
\end{proof}

\begin{mtheorem}{\bf Theorem.}
\label{theorem-primitive-3}
Let $X$ be a primitive del Pezzo threefold with $\rr(X)=3$.
Let $\{\FFF_i\}$ be the set of all pencils as in Lemma \xref{lemma-primitive-r=3}.
Then there are the following possibilities for $\{\FFF_i\}$, where we 
draw the graph for $\{\FFF_i\}$
so that every two elements are connected by an edge if and only if they are conjugate.

\begin{enumerate}
 \item \label{primitive-rho=3-d=2}
{$\dd(X)=2$}

\[
\xygraph{
!{<0cm,0cm>;<2.8cm,0cm>:<0cm,1cm>::}
!{(0,3) }*+{\overset{S-\FFF_1+\FFF_2}\bullet}="S-F1F2"
!{(1,2) }*+{\overset{2S-\FFF_1}\bullet}="2SF1"
!{(-1,2) }*+{\overset{\FFF_2}\bullet}="F2"
!{(1,1) }*+{\underset{S-\FFF_2}\bullet}="2SF2"
!{(-1,1) }*+{\underset{\FFF_1}\bullet}="F1"
!{(0,0) }*+{\underset{S+\FFF_1-\FFF_2}\bullet}="SF1-F2"
"S-F1F2"-"F2"
"S-F1F2"-"2SF1"
"2SF1"-"2SF2"
"F2"-"F1"
"2SF2"-"SF1-F2"
"F1"-"SF1-F2"
} 
\]

\item \label{primitive-rho=3-d=4}
{$\dd(X)=4$} 
\[
\xygraph{
!{<0cm,0cm>;<2.5cm,0cm>:<0cm,1cm>::}
!{(0,0) }*+{\underset{S-\FFF_2}\bullet}="SF2"
!{(1,0) }*+{\underset{S-\FFF_1}\bullet}="SF1"
!{(0,1) }*+{\overset{\FFF_1}\bullet}="F1"
!{(1,1) }*+{\overset{\FFF_2}\bullet}="F2"
"SF2"-"SF1"
"SF2"-"F1"
"F2"-"SF1"
"F2"-"F1"
} 
\]

\item
\label{primitive-rho=3-d=6}
{$\dd(X)=6$ and $X\simeq \PP^1\times \PP^1\times\PP^1$} 
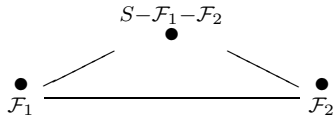
\[
\xygraph{
!{<0cm,0cm>;<2cm,0cm>:<0cm,1cm>::}
!{(0,0) }*+{\underset{\FFF_1}\bullet}="a"
!{(1,1) }*+{\overset{S-\FFF_1-\FFF_2}\bullet}="b"
!{(2,0) }*+{\underset{\FFF_2}\bullet}="c"
"a"-"b" "a"-"c"
"b"-"c"
} 
\]
\end{enumerate}
\end{mtheorem}

\begin{proof}
Let $\FFF_1$ and $\FFF_2$ be two conjugate pencils and let 
$\xi :\hat{X}\to X$ be the corresponding small $\QQ$-factorialization.
Clearly, we have
\[
\FFF_1^2\equiv \FFF_2^2\equiv 0,\ 
\FFF_1\cdot \FFF_2\cdot S=1,
\ S^2\cdot \FFF_1=S^2\cdot \FFF_2=2. 
\]
For any $j$, write $\FFF_j\sim aS+b_1\FFF_1+b_2\FFF_2$, where $a\ge 0$. Then
\begin{equation}
\label{eq-primitive-rk=3}
\begin{array}{ll}
0=\FFF_j^2\cdot S&=a^2d+4a(b_1+b_2)+2b_1b_2,
\\[5pt]
2=\FFF_j\cdot S^2&=ad+2(b_1+b_2),
\end{array}
\end{equation}
where $d:=\dd(X)$. Therefore,
\[
\begin{array}{lll}
b_1+b_2&=&\frac12(2- ad),
\\[5pt]
b_1b_2&=&\frac12 a(ad-4).
\end{array}
\]
Since this system has an integer solution in $b_1$, $b_2$, the 
discriminant 
\[
\frac14(2- ad)^2-2 a(ad-4)=\frac 14\bigl(4- a (8 - d) (a d - 4)\bigr)
\]
must be a square and $ad$ must be even. Assuming $a>0$ (i.e. $\FFF_j\neq \FFF_1$, $\FFF_2$), 
we get $ad=8$, $6$, $4$, or $2$. 
Hence, up to permutation of $b_1$ and $b_2$, there are 
the following solutions with $a>0$:

\begin{enumerate}
 \item[]
$d=1$,\quad $(a, b_1, b_2)=$ 
 $(4, -1, 0)$,
 $(4, 0, -1)$;
\item[]
$d=2$,\quad $(a, b_1, b_2)=$
 $(1, -1, 1)$,
 $(1, 1, -1)$,
$( 2, -1, 0)$,
 $(2, 0, -1)$;
\item[]
$d=4$,\quad $(a, b_1, b_2)=$
$( 1, -1, 0)$,
 $(1, 0, -1)$;
\item[]
$d=6$,\quad $(a, b_1, b_2)=$
 $(1, -1, -1)$.
\end{enumerate}
Note that if $\FFF_j$ and $\FFF_k$ are conjugate, then 
$\FFF_j\cdot \FFF_k\cdot S=1$. 
From this one can see that for each $\FFF_j$ 
there are exactly two divisors in $\{\FFF_i\}$ that conjugate to $\FFF_j$.
Moreover, if $d\neq 1$, then
conjugacy relations are given by graphs in 
Theorem~\ref{theorem-primitive-3} \ref{primitive-rho=3-d=2},
\ref{primitive-rho=3-d=4}, \ref{primitive-rho=3-d=6}.
In the case $d=1$ we get the following (disconnected) graph:
\[
\xygraph{
!{<0cm,0cm>;<3cm,0cm>:<0cm,1cm>::}
!{(0,0) }*+{\overset{\FFF_1}\bullet}="F1"
!{(1,0) }*+{\overset{\FFF_2}\bullet}="F2"
!{<0cm,0cm>;<3cm,0cm>:<0cm,1cm>::}
!{(2,0) }*+{\overset{4S-\FFF_1}\bullet}="4SF1"
!{(3,0) }*+{\overset{4S-\FFF_2}\bullet}="4SF2"
"4SF1"-"4SF2"
"F1"-"F2"
}
\]

Hence there are only two 
extremal $K_{\hat X}$-negative contractions on $\hat X$.
On the other hand, the cone $\overline{\operatorname{NE}}^1(\hat X)$ has 
at least three extremal rays, a contradiction.
\end{proof}

\begin{case}{\bf Remark.}
Let $X$ be a primitive del Pezzo threefold with $\rr(X)=3$ and $\dd(X)=2$ or $4$.
Let $\xi : \hat X\to X$ be a $\QQ$-factorialization. 
Then $\hat X\simeq \PP(\EEE)$, where $\EEE$ is a stable 
rank two vector bundle on $Z=\PP^1\times \PP^1$
with $c_1(\EEE)=0$, $c_2(\EEE)=6-\dd(X)$.

\end{case}

\begin{scase}{\bf Example.}
If $\dd(X)=4$, an example of such $\EEE$ can be obtained as a restriction of 
the null-correlation bundle $\mathscr N$ from $\PP^3$ to $Z$, where $Z\subset\PP^3$
is the Segre embedding. Recall that the null-correlation bundle is defined by 
the exact sequence
\[
0\longrightarrow \OOO_{\PP^3}\longrightarrow
\Omega_{\PP^3}(2) \longrightarrow\mathscr N(1) \longrightarrow 0.
\]
Its projectivization $Y:=\PP(\mathscr N)$ is a Fano fourfold of index $2$ \cite{Szurek1990}. 
This $Y$ has also a structure of $\PP^1$-bundle over a smooth three-dimensional quadric.
Let $\hat X=\PP(\EEE)=\pi^{-1}(Z)$, where $\pi :Y\to \PP^3$ is the natural projection.
Then $\hat X$ is a weak del Pezzo threefold of type \xref{theorem-primitive-3}\ref{primitive-rho=3-d=4}.
\end{scase}
Examples of del Pezzo threefolds of type~\ref{theorem-primitive-3}
\ref{primitive-rho=3-d=2} can be constructed similarly 
by restricting to $Z\subset \PP^3$ rank two stable vector bundles $\mathscr F$ 
with $c_1=0$, $c_2=2$ \cite[\S 9]{Hartshorne1978}. 

Another way to show existence of del Pezzo threefolds of types~\ref{theorem-primitive-3}\ref{primitive-rho=3-d=4}
and \ref{primitive-rho=3-d=2}
is by writing down explicit equations:
\begin{scase}{\bf Example.}
Let $X\subset \PP^5$ is given by the equations
\[
\begin{cases}
x_1x_3-x_2x_4+a_{3,4}x_3x_5 +a_{3,6}x_3x_6+a_{4,5}x_4x_5+a_{4,6}x_4x_6=0
\\
x_1x_5-x_2x_6+b_{3,4}x_3x_5 +b_{3,6}x_3x_6+b_{4,5}x_4x_5+b_{4,6}x_4x_6=0
\end{cases}
\]
where $a_{i,j}$, $b_{i,j}$ are sufficiently general constants. 
Then $X$ is a del Pezzo threefold having exactly $4$ nodes.
By Corollary \ref{corollary-singular-points-V14} $\rr(X)\ge 3$.
On the other hand, by results of \ref{maximal-deg=4} and \S \ref{section-submaximal}
below $\rr(X)=3$. Finally, two quadrics $x_5=x_6=x_1x_3-x_2x_4=0$ and
$x_3=x_4=x_1x_5-x_2x_6=0$ determine two conjugate pencils.
Therefore, $X$ is of type \xref{theorem-primitive-3}\ref{primitive-rho=3-d=4}.
\end{scase}
 \comment{
\begin{scase}{\bf Example.}
Let $X\subset \PP(1,1,1,1,2)$ is given by the equation
\[
x_0^2=q^2-(x_1x_2-x_3x_4)(x_1x_3-x_2x_4),
\]
where $\deg x_0=2$, $\deg x_i=1$ and $q=q(x_1,\dots,x_4)$ is a general 
quadratic form. Then $X$ has exactly $12$ nodes
and does not contain any plane. 
Hence by Table \ref{main-theorem} $X$ is of type~\ref{theorem-primitive-3}\ref{primitive-rho=3-d=2}.
\end{scase}
}

\section{Del Pezzo threefolds with $\rr(X)=2$}
\label{section-Primitive-2}
The results of this section are contained in 
\cite{Jahnke2008}. We give a short self-contained proof for the 
convenience of the reader. 
\begin{case}
Let $X$ be a del Pezzo threefold with $\rr(X)=2$.
There exists the following diagram:
\[
\xymatrix{
&\hat{X}\ar[dl]_{f}\ar[dr]_{\xi}\ar@{-->}[rr]
&&\hat{X}^+\ar[dl]^{\xi^{+}}\ar[dr]^{f^+}
\\
Z&&X&&Z^+
}
\]
where $\xi$, $\xi^+$ are small $\QQ$-factorializations,
$\hat X \dashrightarrow \hat X^+$ is a flop, and 
$f$, $f^+$ are $K$-negative extremal contractions.
We may assume that $\dim Z\ge \dim Z^+$.
Let $S=-\frac12 K_X$ and let $\hat S=h^*S$.
Let $M$ (resp. $M^+$) be the ample generator of $\Pic(Z)$
(resp. $\Pic(Z^+)$). Put $L:=f^*M$ and $L^+:=f^{{+}*}M^+$.
Let $L'$ be the proper transform of $L^+$ on $\hat X$.
If $f$ is birational, then $E\subset \hat X$ 
denotes the $f$-exceptional divisor.
Similarly, if $f^+$ is birational, then 
$E'\subset \hat X$ is the 
proper transform of $f^+$-exceptional divisor.
\end{case}

\begin{scase}{\bf Remark.}
\label{remark-flop}
If in the above notation $\dd(X)\le 2$, then by
\xref{def-Bertini-Geiser} there is a natural (Bertini or Geiser) involution
$\tau: X\to X$. In this case,  $\hat X^+\simeq \hat X$,
$\xi^+=\tau\comp \xi$,   $Z\simeq Z^+$, and $f^+$ has the same type as $f$. 
\end{scase}

The following theorem was proved (in much stronger form) in \cite{Jahnke2008}.
A short simple proof can be found in the arxiv version of the paper.
\begin{mtheorem}{\bf Theorem.}
\label{theorem-primitive-2}
In the above notation 
there are the following possibilities.
\begin{center}
\setlongtables
\begin{longtable}{c|c|c|c|c|c}
\hline
 &$f$ & $f^+$ &$\dd$ & $\Pic(\hat X)$&$\s$
\\[4pt]
\hline
\endfirsthead
 &$f$ & $f^+$ &$\dd$ & $\Pic(\hat X)$&$\s$
\\
\hline
\endhead
\endlastfoot
\endfoot
&&&&&
\\[-7pt]
\numer\label{rkCl=2-P1-bundle-degree=1}
&$\PP^1$-bundle&$\PP^1$-bundle &$1$&$L+L'\sim 6\hat S$& $22$
\\
\numer\label{rkCl=2-P1-bundle-degree=2}
&&&$2$&$L+L'\sim 3\hat S$& $11$
\\
\numer
\label{rkCl=2-P1-bundle-degree=3}
&&&$3$&$L+L'\sim 2\hat S$ &$6$
\\
\numer\label{rkCl=2-P1-bundle-degree=6}
&&&$6$&$L+L'\sim \hat S$&$0$
\\
\hline
&&&&&
\\[-7pt]
\numer\label{rkCl=2-degree=5}
&$\PP^1$-bundle & quadric bundle &$5$&$L+L'\sim \hat S$&$1$
\\
\hline
&&&&&
\\[-7pt]
\numer\label{rkCl=2-deg1-quadric}
&quadric bundle & quadric bundle &$1$&$L+L'\sim 4\hat S$& $\le 22$
\\
\numer\label{rkCl=2-deg2-quadric}
&& &$2$&$L+L'\sim 2\hat S$& $\le 11$
\\
\numer\label{rkCl=2-deg4-quadric}
&& &$4$&$L+L'\sim \hat S$& $\le 3$
\\
\hline
&&&&&
\\[-7pt]
\numer\label{rkCl=2-deg4-birational}
&birational & $\PP^1$-bundle &$4$&$E+L'\sim \hat S$&$3$
\\
\numer\label{rkCl=2-W7}& & 
&$7$&$E+2L'\sim \hat S$&$0$
\\
\hline
&&&&&
\\[-7pt]
\numer\label{rkCl=2-cubic-plane}
&birational & quadric bundle &$3$&$E+L'\sim \hat S$& $4$, $5$, $6$
\\
\hline
&&&&&
\\[-7pt]
\numer\label{rkCl=2-deg1-birational}
&birational & birational &$1$&$E+E'\sim 2\hat S$& $12\le \s\le 22$
\\
\numer\label{rkCl=2-deg2-birational}
&&&$2$&$E+E'\sim \hat S$& $6\le \s\le 11$
\end{longtable}
\end{center}
Here in the 5th column we indicate relations between $L$, $L'$, $E$, and $E'$ in 
$\Pic(\hat X)$.
\end{mtheorem}

\begin{mtheorem}{\bf Corollary.} 
\label{corollary-rkcl=1}
Let $X$ be a del Pezzo threefold with $\rr(X)=1$.
Assume that $X$ is singular.
Then $\dd(X)\le 4$. 
If $\dd(X)= 4$, then every singular point 
$P\in X$ is rs-nondegenerate \textup(see \xref{Definition-Appendix}\textup). 
Moreover, 
$\sum_P \lambda(X,P)\le 2$.
\end{mtheorem}

\begin{proof}
Let $P\in X$ be a general point. Let $\sigma: \tilde X \to X$ be
the blowup of $P$, let $E:=\sigma^{-1}(P)$, and 
let $\tilde S$ be the proper transform of $S=-\frac12K_X$.
Write $-K_{\tilde X}=2\sigma^*S-2E=2\tilde S$.
Since the linear system $|\tilde S|$ is base point free and big, 
$\tilde X$ is a weak del Pezzo threefold 
with at worst factorial terminal singularities, $\rho(\tilde X)=2$, and
$\dd(\tilde X)=\dd(X)-1$.
If $\dd(X)\ge 5$, then by Theorem \ref{theorem-primitive-2} 
we have only one possibility \ref{rkCl=2-deg4-birational}.
But then both $\tilde X$ and $X$ are smooth.
If $\dd(X)=4$, then we have case
\ref{rkCl=2-cubic-plane}.
In this case any singularity $\tilde P\in \tilde X$ is analytically isomorphic to
the hypersurface singularity given by $x_1x_2+x_3^2+x_4^{n}=0$.
Then $\lambda(\tilde X,\tilde P)=\nu(\tilde X,\tilde P)=\lfloor n/2\rfloor$.
The last inequality follows by Proposition \ref{claim-singular-points-V14}.
\end{proof}

\begin{case}
By \cite{Jahnke2008} all the cases in the table do occur\footnote{There is a typographical error in 
\cite[Th. 3.6]{Jahnke2008}: the case $c_2(\mathcal F)=6$ occurs.}.
Below we give explicit examples of some del Pezzo threefolds 
with $\rr(X)=2$.

\subsection*{Case \ref{rkCl=2-P1-bundle-degree=3}.}
$X=X_3\subset \PP^4$ is given by an equation of the form
\[
(x_1x_4-x_2x_3)\ell_1 +(x_2^2-x_1x_3)\ell_2 +(x_3^2-x_2x_4)\ell_3=0,
\]
where $\ell_i(x_1,\dots, x_5)$ are linear forms.

\subsection*{Case \ref{rkCl=2-degree=5}}
(cf. \xref{maximal-deg=5} and \xref{semi-maximal-deg=5}.)
Let $Y$ be the blowup of $\PP^1\times \PP^2$ 
along a smooth curve $C$ of bidegree $(2,1)$. Then $Y$ is a Fano threefold
with $-K_Y^3=38$ and $\rho(Y)=3$ \cite{Mori1981-82}.
Let $S\subset \PP^1\times \PP^2$ be a (unique) effective divisor 
of bidegree $(0,1)$ containing $C$ and let $\tilde S$ be 
the proper transform of $S$ on $Y$. Then $\tilde S\simeq S\simeq \PP^1\times \PP^1$
and $\OOO_{\tilde S}(\tilde S)$ is of type $(-1,-1)$.
Therefore, there exists a contraction 
$\varphi: Y\to X$, where $\varphi(\tilde S)$ is a node.
Here $X$ is a quintic del Pezzo threefold as in \ref{rkCl=2-degree=5}.

\subsection*{Case \ref{rkCl=2-deg2-quadric}.}
$X\subset \PP(1^4,2)$ is given by the equation
\[
x_5^2=(x_1x_2-x_3x_4)^2+ (x_1x_2-x_3x_4) q_1(x_1,\dots,x_4)+q_2(x_1,\dots,x_4)^2,
\]
where $q_1$ and $q_2$ are general quadratic forms.

\subsection*{Case \ref{rkCl=2-deg4-quadric}.}
$X\subset \PP^5$ is given by the equations
\[
x_1x_2+x_3x_4+x_5^2+x_6l_1(x_1,\dots x_6)= x_1x_3+x_6l_2(x_1,\dots x_6)=0,
\]
where $l_i$ are linear forms.
It is easy to see that $X$ contains 
two singular quadrics given by $x_6=x_1=x_3x_4+x_5^2=0$
and $x_6=x_3=x_1x_2+x_5^2=0$. They generate 
two pencils. Hence $X$ is of type \ref{rkCl=2-deg4-quadric}.
For a general choice of $l_i$ the variety $X$ has exactly one node.

\subsection*{Case \ref{rkCl=2-deg4-birational}.}
$X\subset \PP^5$ is given by the equations
\[
x_3x_4-x_5^2+x_6l_1(x_1,\dots x_6)= x_1x_4-x_2x_5+x_6 l_2(x_1,\dots x_6)=0,
\]
where $l_i$ are general linear forms. Then  $X$ contains the plane $\Pi:=\{x_4=x_5=x_6=0\}$.
The singular locus of $X$ consists of three nodes contained in $\Pi$ and given there by 
$x_3=l_1=0$ and $x_2=x_3l_2-x_1l_1=0$.

\subsection*{Case \ref{rkCl=2-cubic-plane}.}
Let $X\subset \PP^4$ be
given by the following equation:
\[
x_1u(x_1,x_2,x_3,x_4,x_5)+x_2v(x_1,x_2,x_3,x_4,x_5)=0,
\]
where $u$ and $v$ are quadratic forms.
This cubic contains the plane $\Pi:=\{x_1=x_2=0\}$ and, 
for general $u$ and $v$, the singular locus consists of 
four nodes. The projection from $\Pi$
gives us a quadric bundle structure on $\hat X$ (which is the
blowup of $\Pi$).

\subsection*{Case \ref{rkCl=2-deg1-birational}.}
$X\subset \PP(1^3,2,3)$ is given by the equation
\[
x_5^2=x_4^3+ x_4^2\phi_2 + x_4\phi_4 + \phi_3^2,
\]
where $\phi_i(x_1,x_2,x_3)$ are general homogeneous forms 
of degree $i$.

\subsection*{Case \ref{rkCl=2-deg2-birational}.}
$X\subset \PP(1^4,2)$ is given by the equation
\[
x_5^2=x_1\phi_3(x_1,\dots,x_4)+q(x_1,\dots,x_4)^2,
\]
where $\phi_3$ and $q$ are general homogeneous forms 
of degree $3$ and $2$, respectively.
\end{case}

\section{Root systems}
\label{section-Roots}
\begin{case}
\label{notation-roots}
Let $X$ be a del Pezzo threefold of degree $d=\dd(X)$.
In this section we study the image of the restriction map 
$\iota^*: \Cl(X) \to \Pic(S)$, where $S\in |{-}\frac12K_X|$
is a smooth member contained in the smooth locus of $X$
and $\iota: S\hookrightarrow X$ is an embedding. 
Define $\Delta$ and $\Delta'$ as in \ref{delta-1}.
If $X$ is imprimitive, we apply construction \eqref{eq-construction} with all corresponding notation.
In the primitive case, to unify notation, we put $\sigma=\operatorname{id}$.

Note that $S$ does not pass through singular points of $X$.
Thus we may identify $S$ and $\hat S=\xi^{-1}(S)$.
Let $\bar S:= \sigma(S)$. Then $\bar S$ is a smooth del Pezzo surface,
$\bar S\in |{-}\frac 12 K_{\bar X}|$ and $\sigma_S: S\to \bar S$ is a blowup of 
$\rr(X)-\rr(\bar X)$
distinct points.
Define $\bar \Delta$ and $\bar \Delta'$ for 
$\bar S$ as in \ref{delta-1}.
\end{case}

\begin{mtheorem}{\bf Theorem.}
\label{cor-2-rkcl-res}
\textup{(i)}
In the above notation 
the image $\iota^*\Cl(X)$ is the orthogonal complement to $\Delta'$.
In particular, 
\begin{equation}
\label{equation-roots-rank-Delta}
\rk \Delta'+\rk \Cl(X)+\dd(X)=10.
\end{equation}

\textup{(ii)}
We have $\Delta'=\sigma_S^* \bar \Delta'$.

\textup{(iii)}
According to
possibilities for $Z$ we have the following cases:
\begin{enumerate}
\renewcommand{\labelenumi}{(\alph{enumi})}
\item 
If $Z$ is a point \textup(i.e. $\rho(\bar X)=1$\textup), 
then $\bar \Delta'=\bar \Delta$. Here $\bar \Delta'$ is of type $E_8$, 
$E_7$, $E_6$, $D_5$, $A_4$, $A_1$ in cases 
$\dd(\bar X)=1$,
$2$, $3$, $4$, $5$, and $8$, respectively.

\item 
If $Z\simeq \PP^2$, then 
$\bar \Delta'=\{\alpha \in \bar \Delta \mid 
\alpha\cdot f^* K_Z= 0\}$. Here $\bar \Delta'$
is of type $A_{m}$ \textup(recall that $\dd(\bar X)=1$, $2$, $3$, $5$, or $6$\textup).

\item 
If $Z\simeq \PP^1$, then 
$\bar \Delta'=\{\alpha \in \bar \Delta \mid 
\alpha\cdot C= 0\}$, where $C$ is a conic on $\bar S$. Here $\bar \Delta'$
is of type $D_m$
\textup(recall that $\dd(\bar X)=1$, $2$, or $4$\textup).\footnote
{Cases (b) and (c) overlap for $X$ with $\dd(X)=5$.}

\item 
If $X\simeq (\PP^1)^3$, then 
$\Delta'$ is the subsystem $A_1$ in $\Delta\simeq A_1\times A_2$.

\item 
If $X$ is of type~\ref{theorem-primitive-3}\xref{primitive-rho=3-d=2} or \xref{primitive-rho=3-d=4},
then $\Delta'$ is of type $A_5$ or $A_3$, respectively. 
\end{enumerate}
\end{mtheorem}

\begin{proof}
\begin{case}\label{roots-proof-primitive}
Assume that $X$ is primitive. Then $\hat X=\bar X$ and $\sigma=\operatorname{id}$.
All the statements are obvious if 
$\rr(X)=1$. We assume that $\rr(X)\ge2$.
Let 
$f: \hat{X}\to Z$ be an extremal $K_{\hat{X}}$-negative contraction.
Let $S\in |{-}\frac12K_{\hat{X}}|$
be a smooth member.
Denote by $\delta: S\to Z$ the restriction of $f$ to $S$.
Since $f:\hat{X}\to Z$ is an extremal contraction, 
the image of $\iota^*: \Pic(\hat{X}) \to \Pic(S)$ is generated by 
$\delta^*\Pic (Z)$ and $-K_{S}=-\frac12K_{\hat{X}}|_{S}$.
Clearly, $f: S\to Z$ is surjective. 
Fix a standard basis in $\Pic(S)$ \cite[ch. 8]{Dolgachev-topics}:
\[
\h,\quad\e_1,\dots,\e_n, 
\]
where $n=9-d$ and
\[
 \h^2=1, \quad \e_i^2=-1,\quad \e_i\cdot \e_j=0
 \quad \text{for $i\neq j$.}
\]
Since $\iota^*\Pic(\hat{X})$ is generated by $\delta^*\Pic (Z)$ and $-K_S$, we have
\[
\Delta'=\{\alpha\in \Delta \mid \alpha \cdot \delta^*\Pic(Z)=0\}. 
\]

\begin{scase}
\textbf{Case $Z\simeq \PP^2$ and $f$ is a $\PP^1$-bundle.}
Then $f: S\to \PP^2$ is the blowup of $n=9-d$ points and 
we can choose the basis $\h$, $\e_1,\dots,\e_n$
so that $\h=f^*\OOO_{\PP^2}(1)$ and $\e_1,\dots,\e_n$ are $f$-exceptional.
In this case, $\iota^*\Pic(\hat{X})$ is generated by $\h$ and $K_S$.
Hence $\Delta'=\{\alpha\in \Delta \mid \alpha \cdot \h=0\}$.
Then $\Delta'=\{\e_i-\e_j \mid i\neq j\}$ 
is a root subsystem of rank $n-1$ generated 
by $\e_1-\e_2$,\dots,$\e_{n-1}-\e_n$. 
Thus $\Delta'$ is of type $A_{n-1}$.
\end{scase}

\begin{scase}\label{roots-quadric-bundle}
\textbf{Case $Z\simeq \PP^1$, i.e. $f$ is a quadric bundle.}
Then $n\ge 4$ and $\delta: S\to \PP^1$ is a conic bundle. Let $C$ be a fiber.
By changing the basis $\h$, $\e_1,\dots,\e_n$
we may assume that $C\sim \h-\e_1$.
Then 
$\Delta'=\{\alpha\in \Delta \mid \alpha\cdot C=0\}$, i.e.
$\Delta'$ consists of the following elements:
\begin{itemize}
\item 
$\e_i-\e_j$,\quad $i,j>1$, $i\neq j$.
\item 
$\pm(\h-\e_1-\e_i-\e_j)$,\quad $i,j>1$, $i\neq j$.
\end{itemize}
Simple roots can be taken as follows:
\[
 \h-\e_1-\e_2-\e_3, \quad 
 \e_2-\e_3,\dots, \e_{n-1}-\e_n.
\]
Hence $\Delta'$ is of type $D_{n-1}$ if $n\ge 5$ and $A_3$ if $n=4$.
\end{scase}

\begin{scase}
\textbf{Case $Z\simeq \PP^1\times \PP^1$ and $f$ is a $\PP^1$-bundle.}
Let $\ell_i:=F_i|_S$.
Then we may assume that $\ell_1\sim \h-\e_1$, $\ell_2\sim \h-\e_2$.
$\Delta'$ consists of the following elements:
\begin{itemize}
\item 
$\e_i-\e_j$,\quad $i,j>2$, $i\neq j$.
\item 
$\pm (\h-\e_1-\e_2-\e_i)$,\quad $i>2$.
\end{itemize}
Simple roots can be taken as follows:
\[
\h-\e_1-\e_2-\e_3,\quad
 \e_3-\e_4,\dots, \e_{n-1}-\e_n.
\]
Thus $\Delta'$ is of type $A_{n-2}$.
\end{scase}
This proves our theorem in the case where $X$ is primitive.
\end{case}

\begin{case}
Now consider the case where $X$ is imprimitive.
Obviously, the statement of (iii) follows from \ref{roots-proof-primitive}.
There is a
birational contraction $\sigma :\hat X\to \bar X$, where $\bar X$ is primitive
and $\sigma$ is a composition of blowups of smooth points.
Let $l:=\rr(X)-\rr(\bar X)$, let $E_1,\dots, E_l$ be $\sigma$-exceptional divisors, and
let $\e_i=E_i\cap S$ for $i=1,\dots, l$.
By the above, the statement of our theorem holds for $\bar X$
with root system $\bar \Delta'\subset \bar \Delta\subset \Pic(\bar S)$.
We have a commutative diagram
\[
\xymatrix{
\Pic(\bar S)\ar@{^(->}[r]^{\sigma^*_S} & \Pic(S)
&\simeq&\Pic(\bar S)\oplus \sum_{i=1}^l \e_i\cdot \ZZ
\\
\Cl(\bar X)\ar@{^(->}[r]^{\sigma^*}\ar[u]^{\bar \iota^*} &\Cl(X)\ar[u]^{\iota^*} 
&\simeq&\Cl(\bar X)\ar@<35pt>[u]^{\bar \iota^*}\oplus \sum_{i=1}^l E_i\cdot \ZZ\ar@{=}@<-35pt>[u]
}
\]
Now it is easy to see that $\iota^*\Cl(X)^\perp\subset \sigma_S^*\Pic(\bar S)$. Therefore,
\[
\sigma_S^*\bar \Delta'\subset \Delta\cap \iota^*\Cl(X)^{\perp} \subset \Delta\cap \sigma_S^*\Pic(\bar S).
\]
On the other hand, $\sigma_S^*\bar \Delta'\supset \Delta\cap \sigma_S^*\Pic(\bar S)$. 
Hence, $\sigma_S^*\bar \Delta'=\Delta\cap \iota^*\Cl(X)^{\perp}$.
This proves (ii). 
As a consequence we have that the left hand side of \eqref{equation-roots-rank-Delta}
is preserved under birational contractions $\sigma$. 
By \ref{roots-proof-primitive} the equality \eqref{equation-roots-rank-Delta}
holds for primitive del Pezzo threefolds. Thus \eqref{equation-roots-rank-Delta}
holds for imprimitive ones as well. This proves (i).
\end{case}
\end{proof}

\section{Del Pezzo threefolds with maximal $\rr(X)$}
\label{section-maximal}
Recall that $\rr(X)+\dd(X)\le 9$ by Corollary \ref{cor-2-rkcl}.
In this section we study del Pezzo threefolds with $\rr(X)+\dd(X)=9$.

We say that points $P_1,\dots, P_n\in \PP^3$ 
are \textit{in general position} if no three of them 
lie on one line and no four of them lie on one plane.
 
\begin{mtheorem}{\bf Theorem.}
\label{main-theorem-maximal}
Let $X$ be a del Pezzo threefold with $\rr(X)+\dd(X)=9$.
Assume that 
$X\not \simeq \PP^1\times\PP^1\times\PP^1$. Then 

\begin{enumerate}
 \item \label{main-theorem-maximal-0}
$X$ can be obtained by applying construction \eqref{eq-construction} to $\PP^3\simeq V_8\subset \PP^9$
where $\sigma$ is the blowup of 
$n:=\rr(X)-1$ points $P_1,\dots, P_n\in V_8$ in general position. 

 \item \label{main-theorem-maximal-1}
Singular points of $X$ are images of proper transforms of
\begin{enumerate}
\renewcommand{\labelenumi}{\rm (\alph{enumi})}
 \item 
lines passing through $P_i$ and $P_j$, $i\neq j$,
 \item 
twisted cubics passing through six distinct points
$P_{i_1},\dots, P_{i_6}$ \textup(see Claim \xref{claim-twisted-cubics} below\textup). 
\end{enumerate}
 \item\label{main-theorem-maximal-3}
If all the singularities of $X$ are nodes, then 
$\s(X)=28$, $16$, $10$, $6$, $3$, $1$ in cases
$\dd(X)=1$, $2$, $3$, $4$, $5$, $6$, respectively. 
 \item\label{main-theorem-maximal-2}
If $\dd(X)\ge 2$, then all the singularities of $X$ are nodes. 
\end{enumerate}
Conversely, assume that $X$ is a del Pezzo threefold 
whose singularities are at worst nodes and assume that 
$\s(X)=28$, $16$, $10$, $6$, $3$, $1$ in cases
$\dd(X)=1$, $2$, $3$, $4$, $5$, $6$, respectively. 
Then $\dd(X)+\rr(X)=9$.
\end{mtheorem}

Note that in the case $\dd(X)=1$ the statement of \ref{main-theorem-maximal-2}
is wrong: one can easily construct $X$ having only 27 
singular points, where one of them is not a node. 

\begin{stheorem}{\bf Corollary.}
Let $X$ be a del Pezzo threefold with $\rr(X)+\dd(X)=9$. 
If $\dd(X)\ge 3$ and $\dd(X)\neq 6$, then $X$ is unique up to isomorphism.
If $\dd(X)=2$ \textup(resp. $\dd(X)=1$\textup), then $X$ belongs to a $3$-dimensional 
\textup(resp. $6$-dimensional\textup) family.
There are exactly two isomorphism classes of Pezzo threefolds with
$\dd(X)=6$, $\rr(X)=3$.
\end{stheorem}

\begin{proof}
\ref{main-theorem-maximal-0}
If $X$ is primitive, then either $X\simeq \PP^1\times\PP^1\times\PP^1$
or $X\simeq \PP^3$ by Theorems 
\ref{th-del-Pezzo-2},
\ref{theorem-primitive-3}, and \ref{theorem-primitive-2}.
Thus we assume that $X$ is imprimitive and $\dd(X)\le 7$.
We use notation of Theorem \ref{structure-theorem-imprimitive}.
Run construction \eqref{eq-construction}
in such a way that $n$ is maximal possible.
On the last step we get a primitive weak del Pezzo threefold
$\bar X$ with $\rho(\bar X)=9-\dd(\bar X)$. Moreover, 
if $\rho(\bar X)=3$, then $n=0$, $\rho(\hat X)=\rr(X)=3$,
and $\dd(X)=6$.
By Theorem \ref{theorem-primitive-3} we have $X\simeq \PP^1\times\PP^1\times\PP^1$.
On the other hand, by Theorem \ref{theorem-primitive-2}
$\rho(\bar X)\neq 2$.
Hence $\rho(\bar X)=1$, $\dd(\bar X)=8$, and then $\bar X\simeq \PP^3$.

It remains to show that the centers $P_1,\dots,P_n$ of the blowup $\hat X\to \bar X\simeq \PP^3$
are in general position. Indeed, if distinct points $P_i$, $P_j$, $P_k$ lie on 
a line $L\subset \PP^3$, then for its proper transform $L'$ on $\hat X$ we have 
$-K_{\hat X}\cdot L'=-K_{\PP^3}\cdot L-3\cdot 2<0$, a contradiction.
Similarly, if four distinct points $P_i$, $P_j$, $P_k$, $P_l$ lie on 
a plane $D\subset \PP^3$, then
then for its proper transform $\hat D$ on $\hat X$ we have
$K_{\hat X}^2\cdot \hat D=K_{\PP^3}^2\cdot D-4\cdot 4=0$.
Hence $\hat D$ is contracted by the anticanonical map, a contradiction.
This proves \ref{main-theorem-maximal-0}.

\ref{main-theorem-maximal-1}
Let $P\in X$ be a singular point. 
Then $\xi^{-1}(P)$ is a curve. Let $\hat C\subset \xi^{-1}(P)$ be a component and
let $\bar C:=\sigma(\hat C)\subset \bar X$. 
There are two members $\hat S',\, \hat S''\in |{-}\frac12K_{\hat X_{0}}|$
such that $C \subsetneq\hat S'\cap\hat S''$. Then $\bar C \subsetneq\bar S'\cap\bar S''$,
where $\bar S'$, $\bar S''\subset \hat X_{0}=\PP^3$ 
are proper transforms of $\hat S'$ and $\hat S''$.
Therefore, $\deg \bar C\le 3$ and $\bar C$ is not a plane cubic. 
If $\deg \bar C=2$, then $\bar C$ is a conic and it must contain 
four distinct points from $P_1,\dots, P_n$. This contradicts 
our assumption that $P_1,\dots, P_n$ are in general position.
Therefore, $\bar C$ is either a line or a twisted cubic. 
This proves \ref{main-theorem-maximal-1}.

\ref{main-theorem-maximal-3} follows by Corollary \ref{corollary-singular-points-V14}.

\ref{main-theorem-maximal-2}
If $\dd(X)\ge 3$, then $X$ is unique up to isomorphism and the statement 
\ref{main-theorem-maximal-2} can be checked directly (see 
\ref{maximal-deg=6}-\ref{Segre-cubic} below).
Let $\dd(X)=2$
the $\xi$-exceptional set 
consists of proper transforms of lines $L_{i,j}$ passing through 
pairs of distinct points $P_i$, $P_j$ and one twisted cubic $C$
passing through $P_1,\dots, P_6$.
Moreover, the lines $L_{i,j}$ meet 
$C$ transversely.
By blowing the points $P_1,\dots, P_6$ up we get these curves disjointed.
Thus $\xi$ is a small resolution whose
exceptional set is a disjointed union of $16$ smooth rational curves.

The last assertion follows by Corollary \ref{corollary-singular-points-V14}.
\end{proof}

\begin{stheorem}{\bf Claim.}
\label{claim-twisted-cubics}
Let $P_1,\dots, P_6\in \PP^3$ be a points in general position.
Then there exists a twisted cubic curve $C=C_3\subset \PP^3$
containing $P_1,\dots, P_6$. Such a curve is unique.
\end{stheorem}
\begin{proof}
It is easy and left to the reader. 
\end{proof}

By Theorem \ref{cor-2-rkcl-res} we have the following.

\begin{stheorem}{\bf Corollary.}
\label{corollary-maximal-roots}
Let $X$ be a del Pezzo threefold with $\rr(X)=9-\dd(X)$
and $\dd(X)\le 5$. Then the image of $\iota^*: \Cl(X)\to \Pic(S)$ 
is a sublattice orthogonal to some root $\alpha\in \Delta$, i.e.
$\Delta'=\{\pm \alpha\}$. Moreover, $\Delta''$ is of type $E_7$, $D_6$, $A_5$, $A_1\times A_3$, $A_2$
in cases $\dd(X)=1$, $2$, $3$, $4$, $5$, respectively.
\end{stheorem}

\begin{stheorem}{\bf Corollary.}
\label{Corollary-max-G}
 Let $X$ be a del Pezzo threefold with $\rr(X)=9-\dd(X)$
and $\dd(X)\le 4$. 
\begin{enumerate}
 \item 
If $\dd(X)\neq 2$, then the image of the natural map $G\to \Aut(\Delta'')$
is contained in the Weyl group $\W(\Delta'')$.
 \item 
If $\dd(X)\le 3$ and $\Bbbk$ is algebraically closed \textup(i.e. we are in the 
geometric case\textup), then the map $G\to \Aut(\Delta'')$ is an embedding.
\end{enumerate}
\end{stheorem}

\begin{proof}
(i) Similar to \cite[Ch. 4, 26.5]{Manin-Cubic-forms-e-II}. If $\dd(X)=1$, then $\Delta''$ is of type $E_7$ and 
$\Aut(\Delta'')=\W(\Delta'')$ \cite{Serre1987}. 
 For $\dd(X)=3$  the group 
$\Aut(\Delta'')$ is a direct product of 
$\W(\Delta'')$ and $\pm \operatorname{id}$.
If the image of $G$ is not contained in $\W(\Delta'')$, then 
the element $\tau:=-\operatorname{id}$ can be expressed as 
$gw$, where $g\in G$ and $w\in \W(\Delta'')$. 
Note that any reflection $s\in \W(\Delta'')$ can be extended to an element 
$\Aut (\iota^*\Cl(X))$.
Hence, the action of $\tau$ can be 
extended to an action on $\iota^*\Cl(X)$ so that $\tau(K_S)=g w(K_S)=K_S$. 
Let $E$ be a plane on $X$ and let $\e$ be the class $\iota^*(E)$. 
Then 
$$\textstyle
\tau(\e)= \tau \bigl(\frac 1d K_S+\e\bigr) -\frac 1d \tau(K_S)=-\bigl(\frac 1d K_S+\e\bigr)-\frac 1d K_S=
-\frac 2d K_S-\e.
$$
In particular, $2/d$ must be integral, a contradiction.

(ii) Let $G_0$ be the kernel of the map $G\to \Aut(\Delta'')$.
Then $G_0$ acts trivially on $\Cl(X)$. In particular, the diagram \eqref{eq-construction}
is $G_0$-equivariant. Thus $G_0$ acts on $\bar X=\PP^3$ so that there are 
$\ge 5$ fixed points in general position, the images of $\sigma$-exceptional divisors.
Then $G_0$ must be trivial.
\end{proof}

\begin{mtheorem}{\bf Theorem.}
\label{main-theorem-maximal-planes}
Let $X$ be a del Pezzo threefold with $\rr(X)+\dd(X)=9$.
Assume that $X\not \simeq \PP^1\times \PP^1\times \PP^1$.
Let $\Pi\subset X$ is a plane, let $\hat \Pi\subset \hat X$
be its proper transform, and let 
$\bar \Pi =\sigma(\hat\Pi)\subset \bar X=\PP^3$.
Then $\bar \Pi$ is of one of the following types:
\begin{enumerate}
\item 
$\bar \Pi$ is one of the points $P_i$, $\hat \Pi$ is $\sigma$-exceptional;
\item 
$\bar \Pi$ is a plane passing through three of the points $P_i$;
\item 
$\bar \Pi$ is quadratic cone passing through six of the points $P_i$
so that one of them is the vertex of the cone;
\item 
\textup(only for $\dd(X)=1$\textup)
$\bar \Pi$ is cubic surface passing through all the points $P_i$
 so that four of them are double points;
\item 
\textup(only for $\dd(X)=1$\textup)
$\bar \Pi$ is quartic surface passing through all the points $P_i$
 so that all of them are double points and
one of them is a triple point.
\end{enumerate}
The number of planes on $X$ 
is given by the following table:
\begin{center}
\begin{tabular}{c||c|c|c|c|c|c|c}
$\dd(X)$&$7$&$6$ &$5$&$4$&$3$&$2$&$1$
\\[3pt]
\hline
&&&&&&
\\[-6pt]
$\pp(X)$ &$1$&$2$ &$4$&$8$&$15$&$32$&$126$
\\[3pt]
\end{tabular}\end{center}
\end{mtheorem}
\begin{proof}
It is easy to see that all the subvarieties $\Pi$ 
described in (i)-(v) are planes. So the number of planes 
is at least the number indicated in the table.
On the other hand, for any plane $\Pi\subset X$, the intersection 
$\Pi \cap S$ is a line whose class in $\Pic (S)$
is orthogonal to the root $\alpha\in \Pic(S)$ (see 
Corollary \ref{corollary-maximal-roots}).
Define 
$$
\mathcal E:=\{ e\in \Pic(S) \mid e^2=K_S\cdot e=-1,\, e\cdot \Delta'=0\}.
$$
Thus the number of planes is at most $|\mathcal E|$.

Let $\h$, $\e_1,\dots,\e_{9-d}$ be a standard basis of $\Pic(S)$.
Since cases $n\le 3$ are trivial, we may assume that $n\ge 4$.
Then the Weil group $\W(\Delta)$ 
transitively acts on $\Delta$ \cite[8.2.14]{Dolgachev-topics} and 
we can take it so that $\alpha=\e_1-\e_2$.
Now it is easy to compute $\mathcal E$ (cf. \cite{Dolgachev-topics}).
For example, for $d=6$ we have 
$\mathcal E=\{\e_3,\, \h-\e_1-\e_2\}$, and
for $d=5$ we have
$\mathcal E=\{\e_4,\e_4, \h-\e_1-\e_2, \h-\e_3-\e_4\}$. 
Other cases are similar. For $d=1$ we also can observe that
$\mathcal E=\Delta''+K_S$ and apply Corollary \ref{corollary-maximal-roots}. 
\end{proof}

Below we describe del Pezzo threefolds $X$ with $\rr(X)+\dd(X)=9$ explicitly 
and give examples.
These threefolds
were studied extensively in classical literature (see, e.g., \cite[ch VIII, \S 2]{Semple-Roth-1985}).
We assume that $X$ is singular (otherwise $X\simeq \PP^3$, $V_7$, 
or $\PP^1\times \PP^1\times \PP^1$).

\begin{case}{\bf Sextic del Pezzo threefold.}\label{maximal-deg=6}
Let $X\subset \PP^2\times \PP^2$ be given by the equation
$x_1y_1+x_2y_2=0$. Then $X$ is a del Pezzo threefold 
with $\dd(X)=6$ and $\rr(X)=3$. The singular locus 
consists of one node.
\end{case}

\begin{case}{\bf Quintic del Pezzo threefold (cf. \cite{Todd1930}).}\label{maximal-deg=5}
Let $X\subset \Gr(2,5)$ be an intersection of three general
Schubert subvarieties of codimension one.
Then $X$ is a del Pezzo threefold 
with $\dd(X)=5$ and $\rr(X)=4$. The singular locus 
consists of three nodes.
\end{case}

\begin{case}{\bf Quartic del Pezzo threefold.}
\label{maximal-deg=4}
Let $X\subset\PP^5$ be an intersection of two quadrics having 
$6$ isolated singular points. Then in some coordinate system 
$X$ can be given by the equations
\begin{equation}
\label{equation-maximal-deg=4}
x_1^2-x_2^2=x_3^2-x_4^2=x_5^2-x_6^2.
\end{equation}
In \cite[ch VIII, 2.31]{Semple-Roth-1985} this variety is called 
the \textit{tetrahedral quartic threefold}. 
By Corollary \ref{corollary-singular-points-V14} $\rr(X)=5$.
The variety $X$ contains $8$ planes 
\[
\Pi_{\epsilon_1,\epsilon_2,\epsilon_3}=\{
x_1+ \epsilon_1 x_2=x_3 +\epsilon_2 x_4=x_5+\epsilon_3 x_6=0\},
\]
where $\epsilon_i=\pm 1$. Clearly,
\[
\dim \Pi_{\epsilon_1,\epsilon_2,\epsilon_3}\cap 
\Pi_{\epsilon_1',\epsilon_2',\epsilon_3'}=
-1+\frac12\sum |\epsilon_i+\epsilon_i'|.
\]
Therefore, for each plane $\Pi=\Pi_{\epsilon_1,\epsilon_2,\epsilon_3}$
there is exactly $3$ planes $\Pi'$ such that $\Pi\cap \Pi'$ is a point
and exactly $3$ planes $\Pi'$ such that $\Pi\cap \Pi'$ is a line.
Note that 
there are two $4$-tuples of planes such that 
planes in each tuple meet each other 
only by subsets of dimension $\le 0$: 
\[
\{\Pi_{+++},\, \Pi_{+--},\, \Pi_{-+-},\, \Pi_{--+}\},\quad 
\{\Pi_{---},\, \Pi_{-++},\, \Pi_{+-+},\, \Pi_{++-}\}. 
\]
The involution 
\[
\tau: (x_1,x_2,x_3,x_4,x_5,x_6) \longmapsto (x_1,-x_2,x_3,-x_4,x_5,-x_6)
\]
interchanges these $4$-tuples. Hence $\tau$ induces 
a birational (cubo-cubic) involution on 
$\PP^3$. In \cite[Ch. XIV, \S 14, P. 301]{Hudson1927}
it is denoted by $T_{\mathrm{tet}}$. 
Note however that $\Cl(X)^{\tau}\not \simeq \ZZ$, i.e. 
$X$ is not $\tau$-minimal. $X$ is minimal with respect to the whole automorphism group.
\end{case}

\begin{case}\label{Segre-cubic}
{\bf Segre cubic.}
If $\dd(X)=3$, then $X$ is can be given by
\begin{equation}\label{equation-Segre-cubic}
X=X_3^{\mathrm s}=\left\{\sum_{i=1}^6 x_i= \sum_{i=1}^6 x_i^3=0\right\}\subset \PP^4\subset \PP^5. 
\end{equation}
This cubic satisfies many
remarkable properties (see 
\cite[ch VIII, 2.32]{Semple-Roth-1985})
and is called
the \emph{Segre cubic}. For example, any cubic hypersurface in
$\PP^4$ has at most ten isolated singular points, this bound is
sharp and achieved exactly for the Segre cubic (up to projective
isomorphism). The symmetric group $\Sym_6$ acts on $X_3^{\mathrm s}$
in the standard way. Moreover,
by Corollaries \ref{corollary-maximal-roots} and 
\ref{Corollary-max-G} we see that 
$\Aut(X_3^{\mathrm s})=\Sym_6$, so the natural map $\Aut(X_3^{\mathrm s})\to \W(\Delta'')$
is an isomorphism. 
 \end{case}

\begin{case}{\bf Quartic double solid.}
\label{maximal-rk-d=2}
Let $X$ be a del Pezzo threefold of degree $2$.
 Let $\phi: X \to \PP^3$ be the half-anticanonical map.
 Then $\phi$ is a double cover whose branch locus $B\subset \PP^3$
 is a quartic having $16$ singular points.
It is well-known that such a quartic must be a Kummer surface,
so the singularities of $B$ and $X$ are at worst nodes 
\cite{Hudson1905}, \cite{Nikulin1975a} (see also \cite{Jessop1916}).
The singular points of $X$ correspond to $15$ lines $L_{ij}$ passing through 
pairs of points $P_i$, $P_j$ and one twisted cubic passing through 
all points $P_1,\dots,P_6$. The threefold $X$ contains $32$ planes 
\cite[ch VIII, 2.33]{Semple-Roth-1985}. For each such a plane $\Pi$
the image $\pi(\Pi)$ is a plane touching $B$ along a conic.

\begin{scase}{\bf Example.}
Let $B\subset \PP^3$ be a surface given by the 
equation $x_0^4+x_1^4+x_2^4+x_3^4-4x_0x_1x_2x_3=0$.
Then the singular locus of $S$ consists of $16$ isolated points which are simple nodes. 
A double cover $X\to \PP^3$ branched 
along $B$ is a del Pezzo threefold with $\dd(X)=2$ and $\rr(X)=7$. 
\end{scase}
\end{case}
 
\begin{case}{\bf Double Veronese cone.}
\label{maximal-rk-d=1}
Recall that $X\simeq X_6\subset \PP(1^3,2,3)$.
The projection from $(0,0,0,0,1)$ indices a double cover 
$X\to \PP(1^3,2)$ with branch divisor 
$B=B_6\subset \PP(1^3,2)$.
Assume for simplicity that the singularities of $X$
are at worst nodes.
Then $B$ is a surface having exactly $28$ points of type $A_1$.
Conversely if $B\subset \PP(1^3,2)$ is a surface of 
degree 6 whose singularities are exactly $28$ points of type $A_1$,
then the double cover of $\PP(1^3,2)$ branched at $B$
is a del Pezzo threefold with $\dd(X)=1$ and $\rr(X)=8$.
We refer to \cite{Dolgachev-Ortland-1989} for more detailed treatment and more references.

 \begin{scase}{\bf Example.}
 Let $C\subset \PP^2$ is given by the equation
 $f=x_1^4+x_2^4+x_3^4$. Then the dual curve $C^*$ 
 is given by 
 $f^*=(x_1^4+x_2^4+x_3^4)^3-27x_1^4x_2^4x_3^4$.
 It is easy to check that the discriminant of the polynomial
 $h(t)=t^3-(x_1^4+x_2^4+x_3^4)t+2x_1^2x_2^2x_3^2$ is equal to $4f^*$.
 The last polynomial defines a surface $B\subset\PP(1^3,2)$ of degree $6$
 having 28 singular points.
 \end{scase}
 \end{case}

\begin{mtheorem}{\bf Corollary.}
Let $X$ be a del Pezzo threefold such that $\dd(X)+\rr(X)=9$ and 
$\dd(X)\neq 5$, $6$, $7$. Then $X$ is a $G$-del Pezzo threefold with 
respect to some \textup(geometric\textup) group $G$.
\end{mtheorem}

\section{Del Pezzo threefolds with $\rr(X)=8-\dd(X)$}
\label{section-submaximal}
Let, as above, $V_6\subset \PP^7$ be a 
smooth del Pezzo threefold of degree $6$ and 
let $f_i: V_6\to \PP^2$, $i=1$, $2$ be $\PP^1$-bundles. 
We say that points $P_1,\dots, P_n\in V_6$ 
are in \textit{ general position} if so are the points
$f_i(P_1),\dots, f_i(P_n)\in \PP^2$ for $i=1$ and $2$.

\begin{mtheorem}{\bf Theorem.}
\label{main-theorem-semimaximal}
Let $X$ be a del Pezzo threefold with $\rr(X)+\dd(X)=8$.
Then 

\begin{enumerate}
 \item \label{main-semitheorem-maximal-0}
$X$ can be obtained 
by applying construction \eqref{eq-construction} to
$V_6\subset \PP^7$
where $\sigma$ is the blowup of 
from $n:=6-\dd(X)$ points $P_1,\dots, P_n\in V_6$ in general position.

 \item \label{main-semitheorem-maximal-1}
Singular points of $X$ are images of proper transforms of
\begin{enumerate}
\renewcommand{\labelenumi}{\rm (\alph{enumi})}

\item 
curves of bidegree $(0,1)$ and $(1,0)$
passing through one of the points $P_i$;
 
\item 
curves of bidegree $(1,1)$
passing through two of the points $P_i$;

\item 
curves of bidegree $(2,2)$ passing through 
four of the points $P_i$;

\item 
\textup(only for $\dd(X)=1$\textup) 
curves of bidegree $(2,3)$ and $(3,2)$ passing through all the points $P_i$.
\end{enumerate}
 \item\label{main-semitheorem-maximal-3}
If all the singularities of $X$ are nodes, then 
$\s(X)=27$, $15$, $9$, $5$, $2$, $0$ in cases
$\dd(X)=1$, $2$, $3$, $4$, $5$, $6$, respectively. 
 \item\label{main-semitheorem-maximal-2}
If $\dd(X)\ge 2$, then all the singularities of $X$ are nodes. 
\end{enumerate}
Conversely, assume that $X$ is a del Pezzo threefold 
whose singularities are at worst nodes and assume that 
$\s(X)=27$, $15$, $9$, $5$, $2$, $0$ in cases
$\dd(X)=1$, $2$, $3$, $4$, $5$, $6$, respectively. 
Then $\dd(X)+\rr(X)=8$.
\end{mtheorem}

\begin{proof}
Run construction \eqref{eq-construction}
so that $n$ is maximal possible.
On the last step we get a primitive weak del Pezzo threefold
$\bar X$ with $\rho(\bar X)=8-\dd(\bar X)<8$. Moreover, 
if $\rho(\bar X)=3$, then $n=0$, $\rho(\hat X)=\rr(X)=3$,
and $\dd(X)=5$. This is impossible by Theorem \ref{theorem-primitive-3}.
Therefore, $\rho(\bar X)=2$ and $\dd(\bar X)=6$.
By Theorem \ref{theorem-primitive-2} we have only one possibility:
$\bar X\simeq V_6$.
\end{proof}

\begin{stheorem}{\bf Corollary.}
Let $X$ be a del Pezzo threefold with $\rr(X)+\dd(X)=8$.
If $\dd(X)\ge 5$, then $X$ is unique up to isomorphism.
There are exactly two isomorphism classes of 
del Pezzo threefolds with $\dd(X)=\rr(X)=4$.
 \end{stheorem}

\begin{proof}
Indeed, in the case $\dd(X)=4$ two non-isomorphic del Pezzo
threefolds $X$ are obtained by blowing up 
a couple of points corresponding to flags 
$(L_1,P_1)$, $(L_2,P_2)\in \operatorname{F}(\PP^2)=V_6$ such that
either $L_1\cap L_2\neq P_i$ or $L_1\cap L_2= P_i$.
\end{proof}

Similar to Theorem \ref{main-theorem-maximal-planes}
one can prove the following.
\begin{mtheorem}{\bf Theorem.}
\label{main-theorem-submaximal-planes}
Let $X$ be a del Pezzo threefold with $\rr(X)+\dd(X)=8$.
Let $\Pi\subset X$ is a plane, let $\hat \Pi\subset \hat X$
be its proper transform, and let 
$\bar \Pi =\sigma(\hat\Pi)\subset \bar X=V_6$.
Then $\bar \Pi$ is of one of the following types:
\begin{enumerate}
\item 
$\bar \Pi$ is one of the points $P_i$, $\hat \Pi$ is $\sigma$-exceptional;
\item 
$f_j(\bar \Pi)$ is a line for $j=1$ or $2$,
and $\bar \Pi$ contains two of the points $P_i$;
\item 
$\bar \Pi$ is an element of $|{-}\frac12 K_{V_6}|$ 
passing through four of the points $P_i$
so that one of them is a double point;
\item 
\textup(only for $\dd(X)=1$\textup)
$f_j(\bar \Pi)$ is a conic for $j=1$ or $2$,
and $\bar \Pi$ contains all the points $P_i$;

\item 
\textup(only for $\dd(X)=1$\textup)
$\bar \Pi$ is an element of $|{-}K_{V_6}|$ 
passing through all of the points $P_i$
so that all of them are double points
and one of them is triple;

\item 
\textup(only for $\dd(X)=1$\textup)
$\bar \Pi$ is an element of $|{-}K_{V_6}-f_j^*\OOO_{\PP^2}(1)|$,
where $j=1$ or $2$,
passing through all of the points $P_i$
so that three of them are double points.
\end{enumerate}
The number of planes on $X$ 
is given by the following table:
\begin{center}
\begin{tabular}{c||c|c|c|c|c|c}
$\dd(X)$&$6$ &$5$&$4$&$3$&$2$&$1$
\\[3pt]
\hline
&&&&&&
\\[-6pt]
$\pp(X)$ &$0$ &$1$&$4$&$9$&$20$&$72$
\\[3pt]
\end{tabular}\end{center}
\end{mtheorem}

\begin{stheorem}{\bf Corollary.}
\label{corollary-submaximal-roots}
Let $X$ be a del Pezzo threefold with $\rr(X)=8-\dd(X)$
and $\dd(X)\le 5$. Then in some standard basis of $\Pic(S)$
the image of $\iota^*: \Cl(X)\to \Pic(S)$ 
is a sublattice orthogonal to roots $\e_1-\e_2, \, \e_2-\e_3\in \Delta$, i.e.
$\Delta'=\{\pm \e_1\mp\e_2,\, \pm\e_2\mp\e_3, \, \pm (\h -\e_1-\e_2-\e_3) \}$. 
Moreover, $\Delta''$ is of type $E_6$, $A_5$, $2A_2$, $2A_1$
in cases $\dd(X)=1$, $2$, $3$, $4$, respectively.
\end{stheorem}
Similar to Corollary \ref{Corollary-max-G} (ii) we have.
\begin{stheorem}{\bf Corollary.}
 Let $X$ be a del Pezzo threefold with $\rr(X)=8-\dd(X)$
and $\dd(X)\le 2$. 
If $\Bbbk$ is algebraically closed \textup(i.e. we are in the 
geometric case\textup), then the map $G\to \Aut(\Delta'')$ is an embedding.
\end{stheorem}

Now we give some examples. 

\begin{case}{\bf Quintic del Pezzo threefold (cf. \cite{Todd1930}).}\label{semi-maximal-deg=5}
Let $X\subset \Gr(2,5)$ be an intersection of two general
Schubert subvarieties of codimension one and one general hyperplane section.
Then $X$ is a del Pezzo threefold 
with $\dd(X)=5$ and $\rr(X)=3$. The singular locus 
consists of two nodes.

\begin{stheorem}{\bf Corollary (cf. \cite{Todd1930}, \cite{Fujita-1986-1}).} 
\label{corollary-deg=5}
Let $X$ be a del Pezzo threefold of degree $5$.
Then the singularities of $X$ are at worst nodes and 
one of the following holds:
\begin{enumerate}
\item
$X\simeq V_5$, a smooth del Pezzo quintic threefold;
 \item
$\s(X)=1$, $\rr(X)=2$, $\pp(X)=0$, and $X$ is of type 
\xref{theorem-primitive-2}, Case \ref{rkCl=2-degree=5};
 \item
$\s(X)=2$, $\rr(X)=3$, $\pp(X)=1$, and $X$ is of type \xref{semi-maximal-deg=5};
 \item
$\s(X)=3$, $\rr(X)=4$, $\pp(X)=4$, and $X$ is of type \xref{maximal-deg=5}.
\end{enumerate}
\end{stheorem}

\begin{proof}
Assertions (iii) or (iv) 
follows by the results of this and previous sections.
If $\rr(X)=2$, then we have case (ii) by Theorem \ref{theorem-primitive-2}.
Finally, if $X$ is factorial, then it is smooth by Corollary \ref{corollary-rkcl=1}.
\end{proof}
\end{case}

\begin{case}{\bf Quartic del Pezzo threefold.}\label{semi-maximal-deg=4}
Let $X\subset \PP^5$ be given by the equations
\[
x_1^2+x_1x_3+x_2x_5=x_1x_3+x_3^2+ x_4x_6=0.
\]
Then $X$ is a del Pezzo threefold of degree $4$ containing exactly $5$ nodes.
By Corollary \ref{corollary-singular-points-V14} $\rr(X)\ge 4$.
On the other hand, $X$ is not of type \ref{maximal-deg=4} because $\s(X)<6$. 
Hence $\rr(X)=4$.
\end{case}

\begin{case}{\bf Cubic hypersurface.}
Let $X\subset \PP^4$ be given by the equation
\[
x_1x_2\ell(x_1\dots,x_5)+(x_3x_4+x_1x_2)x_5=0,
\]
where $\ell$ is a general linear form.
Then $X$ is a cubic del Pezzo threefold with $\s(X)=9$, 
$\rr(X)=5$, and $\pp(X)=9$ (cf. \cite[J14]{Finkelnberg}).
\end{case}

\begin{case}{\bf Quartic double solid.}
Let $Y$ be 
a hypersurface in $\PP^4$ given by
$\{s_1=4s_4-s_2^2=0\}\subset \PP^5$, where $s_k=\sum x_i^k$. 
This famous hypersurface is called \textit{Igusa quartic}. The singular locus of $Y$ consists of 
$15$ lines. Consider a general hyperplane section $B:=Y\cap \PP^3$.
Then $B$ is a quartic having $15$ nodes (cf. \cite{Jessop1916}).
Let $X\to \PP^3$ be a double cover with branch divisor $B$.
Then $X$ is a del Pezzo threefold of degree $2$ with $\s(X)=15$ and
$\rr(X)=6$.
\end{case}

 \section{$G$-del Pezzo threefolds}
\label{section-G-del-Pezzo}
\begin{case}
In this section we prove Theorem \ref{main-theorem}.
We use notation of \ref{notation-roots}.
Furthermore we assume that $X$ is a $G$-del Pezzo threefold.
Thus $\Cl(X)^G\simeq \ZZ$.
By Theorem \ref{theorem-primitive-2} we may assume that $\rr(X)\ge 3$.
\end{case}

\begin{mtheorem}{\bf Lemma.}
In the above notation, if $\dd(X)\ge 5$, then $X\simeq \PP^1\times \PP^1\times \PP^1$. 
\end{mtheorem}
\begin{proof}
Assume that $X\not \simeq \PP^1\times \PP^1\times \PP^1$.
Then $X$ is singular and $\dd(X)=6$ or $5$ 
by Theorems \ref{th-del-Pezzo-2}.

Consider the case $\dd(X)=6$.
Since $\rr(X)\ge 3$, our $X$ is described in \ref{maximal-deg=6}.
Then $X$ contains exactly two planes $\Pi_1$, $\Pi_2$
and the divisor $\Pi_1+\Pi_2$ is $G$-invariant.
Hence $\Pi_1+\Pi_2\sim a S$ for some positive integer $a$.
Comparing degrees we get $2=6a$, a contradiction.

Now let $\dd(X)=5$. 
By Lemma \ref{lemma-Q-factorial} we may assume that $X$ is not factorial.
In this situation, $X$ is imprimitive. 
The same arguments as above show that the number of planes on $X$ in any $G$-orbit
must be divisible by $5$.
This contradicts Corollary \ref{corollary-deg=5}.
\end{proof}

 \begin{case}
From now on we assume that $\dd(X)\le 4$.
By Theorem \ref{theorem-primitive-3} we may assume that 
$X$ is imprimitive.
Let $S\in |{-}\frac12 K_X|$
be a general member.
Let $n:=\rk \Delta=9-\dd(X)$.

\begin{stheorem} {\bf Lemma.}\label{lemma-planes-degree=4}
If in the above notation $\dd(X)\le 4$, then $X$ contains at least two planes 
$\Pi_1$, $\Pi_2$ such that $\dim \Pi_1\cap \Pi_2\le 0$.
\end{stheorem}

\begin{proof}
Since $X$ is imprimitive, it contains at least one plane $\Pi_1$.
Let $\Pi_1,\dots,\Pi_l$ be its orbit.
Since $\Cl(X)^G=\ZZ\cdot S$, 
$k\ge 4$. If $\dim \Pi_i\cap \Pi_j\ge 1$ for all $i$, $j$,
the linear span of $\Pi_1,\dots,\Pi_k$ is three-dimensional
and so $X$ cannot be an intersection of quadrics.
\end{proof}
\end{case}

First we consider the case $\Delta''=\emptyset$.
\begin{mtheorem} {\bf Proposition.}
\label{construction-d=3-empty-delta}
If in the above notation $\Delta''=\emptyset$, then $\dd(X)=3$, $\rr(X)=3$, $\pp(X)=3$, and 
$X$ is a projection of a del Pezzo threefold $Y=Y_4\subset \PP^5$
of type \eqref{rkCl=2-deg4-quadric} from a point. 
 If moreover the singularities of $X$ are at worst nodes, then 
by \cite{Finkelnberg}, $X$ is of type \textup{J11} or \textup{J12}, and $6\le \s(X)\le 7$.
\end{mtheorem}

\begin{proof}
Let $\e_{n-m+1},\dots, \e_n$ correspond to blowups $\sigma$. 
If $m>1$, then $\e_{n-1}-\e_n\in \Delta''$.
Thus, $m=1$, $\dd(\bar X)=\dd(X)+1$, and we may assume that every 
two planes on $X$ meet each other by a subset of dimension $1$. 
Therefore, $\rr(\bar X)=2$ and $\rr(X)=3$.
By Lemma \ref{lemma-planes-degree=4} $\dd(X)\le 3$. 
Therefore, $\dd(\bar X)\le 4$, and $8\ge n\ge 6$.

Consider the case where $f:\bar X\to Z=\PP^2$ is a $\PP^1$-bundle.
Then by Theorem \ref{cor-2-rkcl-res}
we may assume that the vectors $\e_1-\e_2$, \dots, $\e_{n-2}-\e_{n-1}$
form a basis of $\Delta'$. 
We have:
\begin{itemize}
\item[]
$n=6$, $\dd(X)=3$ $\Longrightarrow$ $2\h-\sum \e_i\in \Delta''$,
\item[]
$n=7$, $\dd(X)=2$ $\Longrightarrow$ $2\h+\e_7-\sum \e_i\in \Delta''$,
\item[]
$n=8$, $\dd(X)=1$ $\Longrightarrow$ $3\h-\e_8-\sum \e_i\in \Delta''$.
\end{itemize}
Thus, in all cases we have $\Delta''\neq \emptyset$, a contradiction. 

Consider the case where $f:\bar X\to Z=\PP^1$ is a quadric bundle.
Then $\dd(X)=\dd(\bar X)-1=1$ or $3$ by Theorem \ref{theorem-primitive-2}.
Again by Theorem \ref{cor-2-rkcl-res}
we may assume that vectors 
\[
\h-\e_1-\e_2-\e_3,\quad
\e_2-\e_3, \dots,\e_{n-2}-\e_{n-1}
\]
form a basis of $\Delta'$. 
If
$n=8$, then $3\h-\e_8-\sum \e_i\in \Delta''$, a contradiction.
Therefore, $\dd(X)=3$ and $\bar X$ is of type \eqref{rkCl=2-deg4-quadric}.
Thus $X$ is a cubic in $\PP^4$. Since for any two planes $\Pi_i,\, \Pi_j\subset X$
we have $\dim \Pi_i\cap \Pi_j\ge 1$, all the planes on $X$ are contained in 
one hyperplane. Hence $\pp(X)=3$.
By Proposition \ref{claim-singular-points-V14}\ $\s(X)\le 7-h^{1,2}(\hat X)$.
If the singularities of $X$ are at worst nodes, then 
$X$ is of type J11 or J12 by \cite{Finkelnberg}.
\end{proof}

\begin{scase}{\bf Example.}
\label{example-d=3-empty-delta}
Consider the cubic $X\subset \PP^4$ given by the equation
\[
x_1x_2x_3+x_0\left(\lambda x_0^2+x_1^2+x_2^2+x_3^2+x_4^2\right)=0. 
\]
Then $X$ has $6$ (resp. $7$) nodes if $\lambda\neq 0$ (resp. $\lambda= 0$). 
It is easy to see that $X$ contains at least $3$ planes, so
$\rr(X)\ge 3$. By \cite{Finkelnberg} we have 
$\rr(X)= 3$ and $\pp(X)=3$. The symmetric group $\Sym_3$
acts on $X$ by permutations of $x_1$, $x_2$, $x_3$ so that $X$ is a $G$-Fano
threefold.
\end{scase}

\begin{case}
Now we assume that $\Delta''\neq \emptyset$.
Then $\Delta''$ is a $G$-invariant 
root subsystem in 
$\Delta$. 
By the results of \S \ref{section-maximal} and \S \ref{section-submaximal}
we may assume that $\rr(X)\le 7-\dd(X)$.
Further, by Lemma \ref{lemma-planes-degree=4} $\dd(X)\le 3$. 
\end{case}

\begin{case}
Consider the case $\dd(X)=3$. There are only the following possibilities:

\begin{scase}
$\dd(\bar X)=4$, $\rr(\bar X)=2$, $\bar X$ is of type \eqref{rkCl=2-deg4-quadric}, $\rr(X)=3$.
Then $\Delta'$ is described in 
\ref{roots-quadric-bundle}: it is of type $D_4$ and generated by 
$\h-\e_1-\e_2-\e_3$, \quad 
 $\e_2-\e_3$,\quad $\e_3-\e_4$,\quad $\e_4-\e_5$.
Any root $\alpha\in \Delta$ has the form $\alpha= \pm (\e_i- \e_j)$,
$\pm (\h-\e_i-\e_j-\e_k)$ or $\pm (2\h-\sum \e_i)$ (see e.g., \cite[ch. 4, 3.7]{Manin-Cubic-forms-e-II}).
Since $\iota^*\Cl(X)=\Delta'^{\perp}$, we get $\Delta''=\emptyset$, a contradiction.
\end{scase}

\begin{scase}
$\dd(\bar X)=5$, $\rr(\bar X)=1$, $\bar X=V_5$, $\rr(X)=3$.
 Similarly, 
$\Delta'$ is of type $A_4$ and generated by
$\h-\e_1-\e_2-\e_3$, \quad 
$\e_1-\e_2$,\quad $\e_2-\e_3$,\quad $\e_3-\e_4$.
In this case, $\Delta''=\{\pm (\e_5-\e_6)\}$.
It is easy to see that the group $G$ permutes elements $\e_5, \,\e_6\in \iota^*\Cl(X)$.
But then the class of $\e_5+\e_6$ must be $G$-invariant, so
it is proportional to $-K_S$, a contradiction.
\end{scase}
\begin{scase}
$\dd(\bar X)=5$, $\rr(\bar X)=2$, $\bar X$ is of type 
\eqref{rkCl=2-degree=5}, $\rr(X)=4$. Similarly, 
$\Delta'$ is of type $A_3$ and generated by
$\e_1-\e_2$,\quad $\e_2-\e_3$,\quad $\e_3-\e_4$.
Then $\Delta''=\{\pm (2\h-\sum \e_i),\ \pm (\e_5-\e_6)\}$.
There is a unique element (class of a line on $S$) $\mathbf x\in (\Delta'+\Delta'')^{\perp}$ 
such that $\mathbf x^2=K_X\cdot \mathbf x=-1$: 
\[
\mathbf x=\h-\e_5-\e_6. 
\]
But then $x\in \iota^*\Cl(X)$ and $x$ must be $G$-invariant, a contradiction.
\end{scase}
\end{case}

\begin{case}
Finally we consider cases $\dd(X)\le 2$. 
According to Remark \ref{def-Bertini-Geiser} any del Pezzo threefold with $\dd(X)\le 2$
is automatically $G$-del Pezzo. Thus all the possibilities for $\bar X$ 
with $2\le \dd(\bar X)\le 5$ and $\rr(\bar X)\le 2$ do occur
(recall that $3\le \rr(X)\le 7-\dd(X)$):

\begin{itemize}
\item 
$\bar X=V_5$ $\Longrightarrow$ $\dd(X)\le 2$, $\Delta'\simeq A_4$;
\item 
$\bar X=V_4$ $\Longrightarrow$ $\dd(X)\le 2$, $\Delta'\simeq D_5$;
\item 
$\bar X=V_3$ $\Longrightarrow$ $\dd(X)\le 2$, $\Delta'\simeq E_6$;

\item 
$\bar X$ is of type \eqref{rkCl=2-P1-bundle-degree=2} $\Longrightarrow$ $\dd(X)=1$, $\Delta'\simeq A_6$;

\item 
$\bar X$ is of type \eqref{rkCl=2-P1-bundle-degree=3} $\Longrightarrow$ $\dd(X)\le 2$, $\Delta'\simeq A_5$;

\item 
$\bar X$ is of type \eqref{rkCl=2-degree=5} $\Longrightarrow$ $\dd(X)\le 2$, $\Delta'\simeq A_3$;

\item 
$\bar X$ is of type \eqref{rkCl=2-deg2-quadric} $\Longrightarrow$ $\dd(X)=1$, $\Delta'\simeq D_6$;

\item 
$\bar X$ is of type \eqref{rkCl=2-deg4-quadric} $\Longrightarrow$ $\dd(X)\le 2$, $\Delta'\simeq D_4$.
\end{itemize}
The number of planes can be found by using Lemma \ref{lemma-plane-restriction}
and direct computations.
\end{case}

\begin{scase}{\bf Example.}
Let $X\subset \PP(1^4,2)$ is given by the equation
\[
y^2= x_1x_2x_3x_4+\lambda (x_1^2+x_2^2+x_3^2+x_4^2)^2
\]
where $\lambda$ is a constant. Then $X$ has exactly $12$ nodes and contains $8$ planes. 
By Corollary \ref{corollary-singular-points-V14} $\rr(X)\ge 3$.
Further, 
by our classification $X$ is of type \ref{class-2-d4}.
\end{scase}
More examples of del Pezzo threefolds with $\dd(X)=2$ can be constructed 
similarly by writing down explicit equations (cf. \cite{Jessop1916}).

\section{Appendix: number of singular points of Fano threefolds}
\label{section-appendix}
\begin{case} {\bf Definition.}
\label{Definition-Appendix}
Let $V\ni P$ be a threefold terminal Gorenstein (=isolated cDV) singularity.
We say that $V\ni P$ is \textit{r-nondegenerate} (resolution nondegenerate)
if there is a resolution 
\[
\sigma:V_m \overset{\sigma_m} {\longrightarrow} 
\cdots \overset{\sigma_2} {\longrightarrow} V_1\overset{\sigma_1} {\longrightarrow} V=V_0,
\]
where each $\sigma_i$ is a blowup of a singular \textit{point} $P_{i-1}\in V_{i-1}$.
Such a resolution $\sigma$ is called \textit{standard}. 
In this situation, all varieties $V_i$ also have only isolated cDV singularities.
If furthermore each $\sigma_i$-exceptional divisor $E_i\subset V_i$ is irreducible,
then we say that $V\ni P$ is \textit{rs-nondegenerate}
(strongly resolution nondegenerate).

Denote $\lambda(V,P):=m$ and let $\nu(V,P)$ be the number of 
$\sigma$-exceptional divisors. Thus $\lambda(V,P)\le \nu(V,P)$
and the equality holds if and only if $V\ni P$ is rs-nondegenerate.
\end{case}

\begin{case} {\bf Remark.}
Let $V\ni P$ be a threefold terminal Gorenstein point
and let $\sigma_1: V_1\to V$ be the blowup of $P$.
Since $V\ni P$ is a hypersurface singularity,
we have an (analytic) embedding $V_1\subset \tilde \CC^4$,
where $\tilde\sigma_1:\tilde \CC^4\to \CC^4$ is the blowup of the origin. 
Let $D:=\tilde\sigma_1^{-1}(P)$ be the exceptional divisor.
Then $D\simeq \PP^3$. Since $V\ni P$ is a singularity of multiplicity $2$,
the intersection $D\cap V_1$  is a quadric in $\PP^3$
(possibly reducible or non-reduced).
If $D\cap V_1$ irreducible, then  $V_1$ is either smooth or has (a unique) terminal singularity.
Moreover, the above arguments show that $2\lambda(V,P)\ge \nu(V,P)$.
\end{case}

\begin{mtheorem}{\bf Proposition.}
Let $(V\ni 0)\subset \CC^4$ be a singularity given by $t^2=\phi(x,y,z)$, where 
$\phi=0$ is an equation of a Du Val singularity. 
Then $V\ni 0$ is r-nondegenerate.
Moreover, if $\phi=0$ defines a singularity of type $A_n$, then 
 $V\ni 0$ is rs-nondegenerate.
\end{mtheorem}

\begin{proof}
Direct computation.
\end{proof}

\begin{stheorem}{\bf Corollary.}
Let $X$ be a del Pezzo threefold with $\dd(X)\le 2$.
Assume that the branch divisor $B$ of the double cover 
$\varphi: X\to \PP(1^{3},2)$ \textup(resp. $\varphi: X\to \PP^3$\textup)
has only Du Val singularities 
\textup(see \xref{def-Bertini-Geiser}\textup).
Then the singularities of $X$ are r-nondegenerate. 
If moreover $B$ has only singularities of type $A$, then the singularities of $X$ are rs-nondegenerate. 
\end{stheorem}

\begin{case}
Let $W$ be a smooth projective fourfold and let $V\subset W$
be an effective divisor. Define
\[
\beta(W,V):=c_3(W)\cdot V-c_2(W)\cdot V^2+c_1( W)\cdot V^3- V^4.
\]
If $V$ is smooth then $\beta(W,V)$ coincides with $\deg c_3(V)=\Eu(V)$,
the topological Euler number.
 \end{case} 
\begin{mtheorem}{\bf Lemma.}
In the above notation let $P\in V$ be a singular point,
let $\sigma: \tilde W\to W$ be the blowup of $P$, 
and let $\tilde V\subset \tilde W$ be the proper transform of $V$.
Then $\beta(\tilde W,\tilde V)=\beta(W,V)+4$. 
\end{mtheorem}

\begin{proof}
Direct computation. See the arxiv version of the paper.
\end{proof}

\begin{mtheorem}{\bf Proposition.}
\label{claim-singular-points-V14}
Let $X$ be a Gorenstein Fano threefold whose singularities are 
r-nondegenerate terminal points. 
Assume that 
\begin{enumerate}
 \item[(*)]
$X$ can be embedded to a smooth fourfold
so that a general member $X'\in |X|$ is smooth.
\end{enumerate}

Then 
\begin{multline*}
{\sum_{P\in X}}'\lambda(X,P)\le 
\sum_{P\in X} \bigl(2\lambda(X,P)-\nu(X,P)\bigr)=
\\
= \rr (X)-\rho(X)+h^{1,2}(X')-h^{1,2}(\hat X), 
\end{multline*}
where $\hat X\to X$ is the standard resolution and the first sum runs through 
all rs-nondegenerate points $P\in X$.
\end{mtheorem}
\begin{proof}
Put $\lambda:=\sum_{P\in X} \lambda(X,P)$.
Thus 
\begin{multline*}
2+2\rho(\hat X)-2h^{1,2}(\hat X)=\Eu(\hat X)=\beta(\hat Y,\hat X)=\beta(\hat Y,\hat X)+4\lambda=
\\
=\Eu(X')+4\lambda=2+2\rho(X')-2h^{1,2}(X')+4\lambda.
\end{multline*}
Since $\rho(\hat X)= \rr(X)+\sum \nu(X,P)$, this gives the desired inequality.
\end{proof}

\begin{scase} {\bf Remark.}
The condition (*) is automatically satisfied if $X$ is a del Pezzo threefold 
(see Theorem \ref{th-del-Pezzo-2}).
\end{scase}

\begin{stheorem}{\bf Corollary.}
\label{corollary-singular-points-V14}
In the notation of
\xref{claim-singular-points-V14} assume additionally that 
the singularities are rs-nondegenerate. 
Then 
\[
|\Sing(X)| \le \rr (X)-\rho(X)+h^{1,2}(X')-h^{1,2}(\hat X).
\]
The equality holds, if all the singularities are 
nodes.
\end{stheorem}

\section{Concluding remarks and open questions}
\label{section-open-questions}
We would like to propose the following open questions. 
\begin{case}{}
Give a complete birational classification of del Pezzo threefolds 
over $\CC$. Non-trivial cases only are 
factorial del Pezzo threefolds of degree $\le 3$.
All other cases can be reduced to the above ones by using 
construction \ref{eq-construction} (or birationally trivial).
It is well-known that a three-dimensional cubic hypersurface 
with at worst cDV singularities is rational if and 
only if it is singular \cite{Clemens-Griffiths}. 
A general smooth (and, in some cases, factorial) del Pezzo threefold of degree $\le 2$
is not rational \cite{Artin-Mumford-1972}, \cite{Beauville1977}, \cite{Tjurin1979}.
\end{case}

\begin{case}{}
Give a complete birational classification of del Pezzo threefolds 
over algebraically non-closed fields. Here is one example.

\begin{stheorem} {\bf Theorem.}
Let $X$ be a smooth del Pezzo threefold of degree $5$
over a field $\Bbbk$. Then $X$ is $\Bbbk$-rational.
\end{stheorem}
\begin{proof}
Denote $\bar X:=X\otimes \bar \Bbbk$. 
Let $\Gamma:=\Gamma(X)$ be the Hilbert scheme parameterizing 
the family of lines on $X$. It is known that 
$\bar \Gamma:=\Gamma\otimes \bar \Bbbk\simeq \PP^2_{\bar \Bbbk}$
(see \cite[Prop. 1.6, ch. 3]{Iskovskikh-1980-Anticanonical}, 
\cite{Furushima1989a}). 
Moreover, lines with normal bundle $N_{l/X}\simeq \OOO(-1)\oplus \OOO(1)$
are parametrized by some conic $C\subset \Gamma$ \cite{Furushima1989a}.
The conic $C$ contains a point of degree $\le 2$. Therefore,
there is a line $\ell \subset \Gamma$ defined over $\Bbbk$.
Let $H_{\ell}$ be the union of all lines $L\subset X$ 
whose class is contained in $\ell\subset \Gamma$.
Then $H_{\ell}$ in an element of $|{-}\frac12 K_X|$ defined over $\Bbbk$
\cite[Proof of Prop. 1.6, ch. 3]{Iskovskikh-1980-Anticanonical}.
In particular, $\Pic(X)=\ZZ\cdot \frac12 K_X$ and the linear system
$|{-}\frac12 K_X|$ defines an embedding $X\subset \PP^6_{\Bbbk}$.
A general pencil of hyperplane sections defines a structure of 
del Pezzo fibration of degree $5$ on $X$. 
By \cite[Ch. 4]{Manin-Cubic-forms-e-II} the variety $X$
is $\Bbbk$-rational.
\end{proof}
\end{case} 

\begin{case}{}
Describe automorphism groups of del Pezzo threefolds over an algebraically closed fields.
Which of them are birationally rigid (cf. \cite{cheltsov-2009}, \cite{Cheltsov2010})? 
These questions are very useful for applications to 
the classification of finite subgroups of Cremona group $\Cr_3(\Bbbk)$
\cite{Prokhorov2009e}, \cite{Prokhorov2009d}.
\end{case}

\def\cprime{$'$} \def\mathbb#1{\mathbf#1}
 \def\bblapr{April}\def\mathbb#1{\mathbf#1}


\end{document}